\documentclass{article}
\usepackage{latexsym} 
\usepackage[all]{xy}
\usepackage{amsmath}
\usepackage{amsfonts}

\def\abstractname{R\'esum\'e}
\newtheorem{thm}{Th\'eor\`eme}[section]
\newtheorem{prop}[thm]{Proposition}
\newtheorem{lem}[thm]{Lemme}
\newtheorem{df}[thm]{D\'efinition}

\newtheorem{cor}[thm]{Corollaire}

\begin{document}

\title{Notes sur la $G$-th\'eorie rationnelle des champs de Deligne-Mumford}
\author{B. Toen\footnote{Max Planck Institut fur Mathematik, Vivatsgasse 7, 53111 Bonn, Germany \newline
 \hspace*{5mm} toen@mpim-bonn.mpg.de}}

\maketitle

\begin{abstract}
Bas\'e sur les m\'ethodes utilis\'ees pour d\'emontrer la formule de Riemann-Roch pour les champs alg\'ebriques (\cite{t1,t2}), 
ce travail donne une description des spectres de $G$-th\'eorie rationnelle des champs de Deligne-Mumford sur
des bases quelconques. Nous utiliserons ces r\'esultats pour \'etudier la $K$-th\'eorie \'equivariante, ainsi que pour d\'efinir de
nouvelles filtrations sur la $K$-th\'eorie des champs alg\'ebriques, permettant potentiellement de d\'evelopper un formalisme de
Riemann-Roch analogue \`a \cite{so}.

\end{abstract}

\def\abstractname{Abstract}

\begin{abstract}
Based on the methods used to prove the Riemann-Roch formula for stacks (\cite{t1,t2}), 
this paper contains a description of the rationnal $G$-theory spectrum of Deligne-Mumford stacks over
general bases. We will use these results to study equivariant $K$-theory, and also to define new filtrations on 
$K$-theory of algebraic stacks, wich allows potentialy to developp a Riemann-Roch formalism analog to \cite{so}.

\end{abstract}

\textsf{Mots cl\'es:} $G$-th\'eorie, champs alg\'ebriques.

\newpage

\tableofcontents

\newpage

\begin{section}{Introduction}

Dans \cite{t1,t2} on a montr\'e comment utiliser la notion de 'cohomologie \`a coefficients dans les
repr\'esentations', afin de d\'emontrer des formules de Riemann-Roch pour les champs alg\'ebriques. L'introduction de cette nouvelle
th\'eorie cohomologique \'etait directement inspir\'ee des r\'esultats d\'ecrivant les spectres de $G$-th\'eorie des 
champs alg\'ebriques (\cite[Thm. $3.15$]{t1}, \cite[Thm. $2.15$]{t2}). Nous donnons ici la version g\'en\'erale de ces descriptions, 
pour le cas des champs de Deligne-Mumford. \\

Nous commencerons dans le premier paragraphe par donner les d\'efinitions n\'ecessaires pour \'enoncer notre th\'eor\`eme. Il s'agit
essentiellement de g\'en\'eralisation des id\'ees utilis\'ees dans \cite{t1}. La petite diff\'erence est que nous utiliserons
le 'champ des sous-groupes cycliques' plutot que le 'champ des ramifications' qui est mal adpat\'e lorsque le champ
ne contient pas assez de racines de l'unit\'e.

La preuve du th\'eor\`eme sera donn\'ee dans le second paragraphe. Ici aussi les id\'ees proviennent directement de \cite{t1}. En
particulier nous utiliserons les th\'eor\`emes de descente, ainsi que la notion de quasi-enveloppe de Chow. 
Nous passerons assez vite sur les d\'etails de la d\'emonstration que le lecteur peut trouver dans \cite{t1,t2}. 

Enfin, nous donnerons deux exemples d'application de ces r\'esultats. 
Tout d'abord nous montrerons comment on peut r\'epondre \`a des questions pos\'ees par A. Vistoli dans \cite{v2}, concernant les
groupes de $G$-th\'eorie \'equivariants, ainsi qu'\`a certaines de leurs g\'en\'eralisations au cas des champs alg\'ebriques.
Enfin, nous d\'efinirons une $\gamma$-filtration sur la $K$-th\'eorie
des champs alg\'ebriques, qui permet hypoth\'etiquement de construire un formalisme de Riemann-Roch analogue \`a celui d\'ecrit
dans \cite{so}. On esp\`ere ainsi g\'en\'eraliser le th\'eor\`eme de Grothendieck-Riemann-Roch de \cite{t1} au cas des champs alg\'ebriques
sur des bases g\'en\'erales. 

Nous terminerons par un bref apper\c{c}u des g\'en\'eralisations \'eventuelles au cas des champs d'Artin. \\

Comme son titre l'indique, le pr\'esent travail est \'ecrit sous forme de notes. Cela implique que nous ne pr\'etendons pas \`a \^etre
exhaustifs, que se soit pour les r\'esultats \'enonc\'es, ou encore pour les arguements qui leur servent de preuve. \\

Ce travail \`a \'et\'e accompli pour les besoins d'un \'expos\'e fait \`a l'institut Max Planck en Octobre $1999$, que je tiens
\`a remercier pour son acceuil. 
Il a aussi fait l'objet d'un expos\'e \`a la conf\'erence de Bologne sur la th\'eorie des
intersections en d\'ecembre $1999$, dont je remercie les organisateurs pour m'avoir permit de m'exprimer sur le sujet. \\

Le pr\'esent travail est la continuation directe de mon travail de th\`ese, dirig\'e par Joseph Tapia et Carlos Simpson, que je 
tiens \`a remercier de nouveau, et sans qui il n'aurait certainement jamais vu le jour.

Je remercie aussi A. Vistoli et G. Vezzosi pour m'avoir signaler la formule \ref{c6'}.

\newpage

\textit{Notations et conventions:} \\

Par la suite nous utiliserons l'expression 'champ alg\'ebrique' pour signifier 'champ alg\'ebrique de Deligne-Mumford
s\'epar\'e et de type fini sur un sch\'ema de base noeth\'erien r\'egulier $S$'. 
Ils seront g\'en\'eralement not\'e par les lettres $F$, $F'$ \dots 
Nous dirons aussi 'sch\'ema' (resp. 'espace alg\'ebrique') 
pour signifier 'sch\'ema (resp. 'espace alg\'ebrique') s\'epar\'e de type fini sur $S$'. Nous adopterons aussi l'abus de language qui consiste
\`a dire qu'un champ alg\'ebrique est quasi-projectif (resp. projectif), si son espace de modules est quasi-projectif (resp. projectif) sur
$S$. \\

Le petit site \'etale d'un champ alg\'ebrique $F$ sera not\'e $F_{et}$. Nous utiliserons 'morphismes de champs'
pour signifier '$1$-morphismes de champs'. La plupart du temps ils seront m\^eme consid\'er\'es comme morphismes dans
la cat\'egorie homotopique des champs. \\

Pour un site $C$, et un pr\'efaisceau en spectres $H$ sur $C$, nous dirons que $H$ est flasque, si pour tout objet $X$ de $C$, 
et tout morphisme
couvrant $U \longrightarrow X$, le morphisme naturel
$$H(X) \longrightarrow \check{H}(U/X,H):=Holim_{[m]\in \Delta^{o}}H(\underbrace{U\times_{X}U \dots \times_{X}U}_{m+1 \; fois})$$
est une \'equivalence faible de spectres. Le th\'eor\`eme de descente cohomologique affirme que tout pr\'efaisceau en spectres
fibrant (pour la structure de cat\'egorie de mod\`eles ferm\'ee de \cite[Thm. $2.53$]{ja}) est flasque. \\

Pour tout pr\'efaisceau en spectres $H$, nous noterons $\mathbb{H}(C,H)$ le spectre des sections globales d'un mod\`ele fibrant pour
$H$. C'est un objet qui est bien d\'efini dans la cat\'egorie homotopique des pr\'efaisceaux en spectres. De fa\c{c}on g\'en\'erale, 
tous les pr\'efaisceaux en spectres seront consid\'er\'es comme des objets de la cat\'egorie homotopique. Ainsi, lorsque nous
parlerons de pr\'efaisceaux en spectres en anneaux il s'agira d'objet en anneaux dans la cat\'egorie homotopique. De m\^eme, lorsque nous
parlerons d'isomorphismes de spectres,  il s'agira d'isomorphismes dans la cat\'egorie homotopique. \\

Pour tout spectre $H$, nous noterons comme il en est l'habitude $H_{\mathbb{Q}}$ le spectre rationnel associ\'e. Plus g\'en\'eralement, 
si $\Lambda$ est un $\mathbb{Q}$-espace vectoriel, on notera $H_{\Lambda}:=H\wedge K(\Lambda,0)$. Si $H$ est un spectre en 
anneaux, et $\Lambda$ une $\mathbb{Q}$-alg\`ebre, $H_{\Lambda}$ est naturellement un spectre  en anneaux. \\

Pour tout champ alg\'ebrique $F$, on notera $\mathbf{G}(F)$ son spectre de $G$-th\'eorie, et 
$\underline{\mathbf{G}}(F):=\mathbb{H}(F_{et},\underline{G}_{\mathbb{Q}})$ son spectre de $G$-cohomologie
rationnelle (\cite[Def. $2.1$]{t2}). Nous noterons aussi $\mathbf{K}(F)$ son spectre de $K$-th\'eorie des complexes parfaits sur $F$ (\cite{th}).
De m\^eme, nous noterons $\underline{\mathbf{K}}(F):=\mathbb{H}(F_{et},\underline{K}_{\mathbb{Q}})$ son spectre de $K$-cohomologie rationnelle. \\

Enfin, nous avertissons le lecteur que nous ne distinguerons pas les foncteurs \`a isomorphismes pr\`es et les vrais foncteurs. Ainsi, 
la correspondance qui \`a un sch\'ema $X$ associe sa cat\'egorie des faisceaux quasi-coh\'erents $Qcoh(X)$ sera consid\'er\'ee comme
un foncteur $Qcoh : Sch/S \longrightarrow Cat$. Le lecteur soucieux de garder cette distingtion utilisera les principes de strictification des pseudo-foncteurs. 

\end{section}

\newpage

\begin{section}{Morphismes de diagonalisation et d'induction}

Si $H \longrightarrow X$ est un sch\'ema en groupes fini et \'etale, nous appellerons
sous-groupe cyclique de $H$, tout sous-sch\'ema en groupes de $H$ qui localement pour la topologie \'etale sur $X$ est isomorphe
\`a $X\times c$, o\`u $c$ est un groupe discret cyclique. \\

Soit $F$ un champ alg\'ebrique. On peut d\'efinir un nouveau champ $C_{F}$ de la fa\c{c}on suivante. Pour un sch\'ema $X$, les objets
de $C_{F}(X)$ sont les couples $(s,c)$, o\`u $s$ est un objet de $F(X)$, et $c$ est un sous-groupe cyclique du sch\'ema en groupes
des automorphismes de $s$ sur $X$. Un morphisme entre $(s,c)$ et $(s',c')$ dans $C_{F}(X)$ est la donn\'ee d'un morphisme
$u : s \rightarrow s'$ dans $F(X)$ tel que $u.c.u^{-1}=c'$. Le lecteur v\'erifiera que ceci d\'efinit bien 
un champ. 

\begin{df}\label{d1}
Soit $F$ un champ alg\'ebrique. Le champ des sous-groupes cycliques mod\'er\'es $C_{F}^{t}$, est le sous-champ  
de $C_{F}$ form\'e des objets $(s,c) \in C_{F}(X)$, o\`u l'ordre de $c$ est premier aux caract\'eristiques de $X$.
\end{df}

On dispose des morphismes $\xymatrix{C_{F}^{t} \ar[r] &  C_{F} \ar[r] &  F}$, o\`u le second oublie le sous-groupe cyclique,
$(s,c) \mapsto s$. Il est facile
de v\'erifier que ces deux morphismes sont repr\'esentables. De plus, le premier est une immersion ouverte, et le second est 
fini et non-ramifi\'e. 

On dispose aussi de la description locale suivante (dont on peut d\'eduire les assertions pr\'ec\'edentes). Notons $p : F \longrightarrow M$
la projection de $F$ sur son espace de modules. On sait que localement sur $M_{et}$, le champ $F$ est \'equivalent \`a un
$[X/H]$, o\`u $X$ est un sch\'ema et $H$ est un groupe fini op\'erant sur $X$. Comme on ne s'int\'eresse qu'\`a la structure locale
on supposera donc que $F=[X/H]$. Notons alors $c(H)$ un ensemble de repr\'esentants de
l'ensemble des classes de conjugaisons de sous-groupes cycliques de
$H$. Pour $c \in c(H)$, notons $X^{c}$ le sous-sch\'ema ferm\'e des points fixes de l'action de $c$ sur $X$, et $U^{c}$ l'ouvert
de $X^{c}$ o\`u l'ordre de $c$ est inversible. Notons aussi pour $c \in c(H)$, $N_{c}$ le normalisateur de $c$ dans $H$, qui
agit naturellement par restriction sur $X^{c}$ et $U^{c}$. On peut alors v\'erifier que les champs $C_{F}$ et $C_{F}^{t}$ sont
donn\'es par
$$C_{F}\simeq \coprod_{c \in c(H)}[X^{c}/N_{c}] \qquad C_{F}^{t}\simeq \coprod_{c \in c(H)}[U^{c}/N_{c}].$$
De part le fait que $C_{F}$ classifie les sous-groupes cycliques des automorphismes de $F$, il existe un champ en groupes universel
$$\mathcal{C}_{F} \longrightarrow C_{F},$$
dont la fibre au-dessus d'un sch\'ema $X \longrightarrow C_{F}$ correspondant \`a $(s,c) \in C_{F}(X)$ est le sch\'ema en groupes
$c \longrightarrow X$. Ce sch\'ema en groupes poss\`ede un faisceau des caract\`eres
$$\begin{array}{cccc} 
\mathcal{X} : & Sch/C_{F} & \longrightarrow & Gp \\
& ((s,c) : X \rightarrow C_{F}) & \mapsto & Hom_{gp/X}(c,\mathbb{G}_{m})
\end{array}$$
Par restriction, on obtient un faisceau en groupes ab\'eliens sur $C_{F}^{t}$, encore not\'e $\mathcal{X}$. Comme la restriction
de $\mathcal{C}_{F}$ sur $C_{F}^{t}$ est un sch\'ema en groupes fini et de type multiplicatif (\cite{sga}), 
le faisceau $\mathcal{X}$ est un faisceau
localement constant pour la topologie \'etale. De plus, il est localement isomorphe au faisceau des racines
de l'unit\'e $\mu_{m}$, pour un certain $m$. 

Pour tout champ alg\'ebrique $F$, on dispose du faisceau associ\'e au pr\'efaisceau en $\mathbb{Q}$-alg\`ebres de groupes
$$\begin{array}{ccc} 
Sch/C^{t}_{F} & \longrightarrow & \mathbb{Q}-Alg \\
((s,c) : X \rightarrow C_{F}^{t}) & \mapsto & \mathbb{Q}[Hom_{gp/X}(c,\mathbb{G}_{m})]
\end{array}$$
Ce faisceau sera not\'e $\mathbb{Q}[\mathcal{X}]$. 
Localement sur $(C_{F}^{t})_{et}$, c'est un faisceau constant isomorphe \`a $\mathbb{Q}[\mathbb{Z}/m]\simeq \frac{\mathbb{Q}[T]}{T^{m}-1}$.
Ainsi, en choisissant $\zeta_{m} \in \mathbb{C}$ une racine primitive de l'unit\'e, il existe localement des quotients
$\mathbb{Q}[\mathcal{X}] \longrightarrow \mathbb{Q}(\zeta_{m})$. On peut v\'erifier que les noyaux de ces morphismes locaux
sont ind\'ependants du choix de $\zeta_{m}$, et aussi du choix de l'isomorphisme $\mathbb{Q}[\mathcal{X}]\simeq \frac{\mathbb{Q}[T]}{T^{m}-1}$.
Ils se recollent donc en un id\'eal $\mathcal{I} \hookrightarrow \mathbb{Q}[\mathcal{X}]$, sur $C_{F}^{t}$.  

\begin{df}\label{d2}
Pour tout champ alg\'ebrique $F$ on notera
$$\Lambda_{F}:=\frac{\mathbb{Q}[\mathcal{X}]}{\mathcal{I}}$$
le faisceau en $\mathbb{Q}$-alg\`ebres quotient. 
\end{df}

Pour en champ alg\'ebrique $F$, le champ $C_{F}^{t}$, muni de son faisceau de $\mathbb{Q}$-alg\`ebres $\Lambda_{F}$,
est 'fonctoriel en $F$'. Ceci signifie que tout morphisme de champs alg\'ebriques $f : F \longrightarrow F'$ 
induit un morphisme naturel $Cf : C_{F}^{t} \longrightarrow C_{F'}^{t}$. De plus, la restriction des caract\`eres permet de d\'efinir
un morphisme naturel de faisceaux de $\mathbb{Q}$-alg\`ebres sur $C_{F}^{t}$
$$Res : Cf^{*}\Lambda_{F'} \longrightarrow \Lambda_{F}$$
En utilisant le morphisme d'induction des caract\`eres, on peut aussi d\'efinir un morphisme de faisceaux en $\mathbb{Q}$-espaces 
vectoriels
$$Ind : \Lambda_{F} \longrightarrow Cf^{*}\Lambda_{F'}.$$
Localement pour la topologie \'etale sur $C_{F}^{t}$, ces morphismes sont isomorphes \`a
$$Res : \mathbb{Q}(\zeta_{m}) \hookrightarrow \mathbb{Q}(\zeta_{m.p})$$
$$Ind=Tr : \mathbb{Q}(\zeta_{m.p}) \longrightarrow \mathbb{Q}(\zeta_{m}),$$
o\`u $Tr$ est le morphisme trace pour l'extension de corps $\mathbb{Q}(\zeta_{m}) \hookrightarrow \mathbb{Q}(\zeta_{m.p})$. \\

Pour tout faisceau en $\mathbb{Q}$-espaces vectoriels $V$ sur un site $C$, on peut construire le pr\'efaisceau en spectres
$K(V,0)$. Vu comme objet de la cat\'egorie homotopique des pr\'efaisceaux en spectres sur $C$, il est caract\'eris\'e par le
fait que $\pi_{m}(K(V,0))=0$ si $m\neq 0$ et $\pi_{0}(K(V,0))\simeq V$. Remarquons que $K(V,0)$ est aussi un pr\'efaisceau en spectres
de Moore (i.e. ses faisceaux d'homologie sont $\underline{H}_{m}(K(V,0))\simeq 0$ si $m\neq 0$ et $\underline{H}_{0}(K(V,0))=V$).  

Si de plus, $V$ est muni d'une structure de faisceau en $\mathbb{Q}$-alg\`ebres, alors $K(V,0)$ devient naturellement un 
pr\'efaisceau en spectres en anneaux.

Soit $H$ un pr\'efaisceau en spectres en anneaux sur $C$. Alors, $H\wedge K(V,0)$ est naturellement un pr\'efaisceau en spectres
en anneaux. De m\^eme, si $K$ est un pr\'efaisceau en spectres qui est un module sur $H$, alors
$K\wedge K(V,0)$ est un module sur $H\wedge K(V,0)$. \\

On peut donner une description explicite et fonctorielle de $H\wedge K(V,0)$ de la fa\c{c}on suivante. \\

Soit $H$ un $\Omega$-spectre, et $A$ un anneau. Notons $e_{n} : H_{n} \longrightarrow \Omega H_{n+1}$ les
morphismes de transitions (qui sont des \'equivalences faibles). On d\'efinit un nouveau spectre par
$$H\otimes A : n \mapsto H_{n}\otimes A,$$ 
o\`u pour un ensemble simplicial $X$, $X\otimes A$ est le $A$-module simplicial $[p] \mapsto X_{p}\otimes A$, obtenu par 
tensorisation. Pour d\'efinir les morphismes de transitions $H_{n}\otimes A \longrightarrow \Omega (H_{n+1}\otimes A)$ il nous suffit 
de d\'efinir un morphisme naturel en $X$
$$(\Omega X)\otimes A \longrightarrow \Omega (X\otimes A).$$
Or, le morphisme naturel $X \longrightarrow X\otimes A$ induit $\Omega X \longrightarrow \Omega (X\otimes A)$, qui 
comme $\Omega (X\otimes A)$ est naturellement un $A$-module simplicial, s'\'etend un le morphisme
cherch\'e. 

On peut v\'erifier que si $A$ est une $\mathbb{Q}$-alg\`ebre alors $H\otimes A$ est naturellement \'equivalent \`a $H\wedge K(A,0)$. Si $H$ poss\`ede
un produit $H\wedge H \longrightarrow H$, on laisse le soin au lecteur de construire naturellement un produit sur $H\otimes A$. \\ 

Appliquons ces derni\`eres remarques \`a $C=(C_{F}^{t})_{et}$ le petit site \'etale du champ $C_{F}^{t}$, 
$V=\Lambda$, $H=\underline{\mathbf{K}}$ le pr\'efaisceau en spectres de $K$-th\'eorie, et $K=\underline{\mathbf{G}}$ le pr\'efaisceau en spectres de $G$-th\'eorie. 
On dispose ainsi du pr\'efaisceau en spectres en anneaux, $\underline{\mathbf{K}}_{\Lambda}:=\underline{\mathbf{K}}\wedge K(\Lambda,0)$, et 
d'un module sur ce dernier $\underline{\mathbf{G}}_{\Lambda}:=\underline{\mathbf{G}}\wedge K(\Lambda,0)$. \\

A l'aide de ces notations, rappelons les d\'efinitions des spectres de $K$-th\'eorie \`a coefficients dans les repr\'esentations. La d\'efinition que
nous donnons ici a subi une l\'eg\`ere modification par rapport \`a \cite[Def. $2.13$]{t2}. 

\begin{df}\label{d3}
Soit $F$ un champ alg\'ebrique. 
\begin{enumerate}

\item Le spectre de $K$-th\'eorie \`a coefficients dans les repr\'esentations de $F$ est d\'efini par
$$\underline{\mathbf{K}}^{\chi}(F):=\mathbb{H}((C_{F}^{t})_{et},\underline{\mathbf{K}}_{\Lambda}).$$

\item Le spectre de $G$-th\'eorie \`a coefficients dans les repr\'esentations de $F$ est d\'efini par
$$\underline{\mathbf{G}}^{\chi}(F):=\mathbb{H}((C_{F}^{t})_{et},\underline{\mathbf{G}}_{\Lambda}).$$
\end{enumerate}
\end{df}

Les $m$-\`eme groupes d'homotopie de ces spectres seront not\'es $\underline{\mathbf{K}}^{\chi}_{m}(F)$ et
$\underline{\mathbf{G}}^{\chi}_{m}(F)$. \\

En utilisant les morphismes de restriction et d'induction d\'efinis pr\'ec\'edemment, on laisse le soin au lecteur de v\'erifier les propri\'et\'es suivantes. \\

\begin{itemize}

\item La correspondance $F \mapsto \underline{\mathbf{K}}^{\chi}(F)$ est un foncteur contravariant de la cat\'egorie homotopique des champs alg\'ebriques vers la cat\'egorie homotopique des spectres en anneaux.

De m\^eme, $F \mapsto \underline{\mathbf{G}}^{\chi}(F)$ est un foncteur covariant de la cat\'egorie homotopique des champs alg\'ebriques et morphismes propres vers celle des spectres. C'est aussi un foncteur contravariant pour les morphismes de $Tor$-dimension finie.

\item Pour tout champ alg\`ebrique $F$, $\underline{\mathbf{G}}^{\chi}(F)$ est un module sur $\underline{\mathbf{K}}^{\chi}(F)$.

De plus, pour tout morphisme propre $p : F \longrightarrow F'$, le diagramme suivant commute 
$$\xymatrix{
\underline{\mathbf{G}}^{\chi}(F) \wedge \underline{\mathbf{K}}^{\chi}(F') \ar[r]^-{p_{*}\wedge Id} \ar[d]_-{Id\wedge p^{*}} & 
\underline{\mathbf{G}}^{\chi}(F') \wedge \underline{\mathbf{K}}^{\chi}(F') \ar[d]^-{-\otimes -} \\
\underline{\mathbf{G}}^{\chi}(F) \wedge \underline{\mathbf{K}}^{\chi}(F) \ar[r]_-{p_{*}(-\otimes -)} &
\underline{\mathbf{G}}^{\chi}(F') \wedge \underline{\mathbf{K}}^{\chi}(F')}$$

\item Pour tour diagramme cart\'esien
$$\xymatrix{ F' \ar[r]^{q} \ar[d]_{g} & G' \ar[d]^{f} \\
F \ar[r]_{p} & G}$$
avec $p$ propre et $f$ plat, le diagramme suivant commute
$$\xymatrix{
\underline{\mathbf{G}}^{\chi}(F) \ar[r]^-{p_{*}} \ar[d]_-{g^{*}} & \underline{\mathbf{G}}^{\chi}(G) \ar[d]^-{f^{*}} \\
\underline{\mathbf{G}}^{\chi}(F') \ar[r]^-{q_{*}} & \underline{\mathbf{G}}^{\chi}(G) }$$

\item Pour tout immersion ferm\'ee $j : F' \hookrightarrow F$, de compl\'ementaire $i : U=F-F' \hookrightarrow F$, on dispose d'une 
suite exacte de fibration
$$\xymatrix{ \underline{\mathbf{G}}^{\chi}(F') \ar[r]^-{j_{*}} & \underline{\mathbf{G}}^{\chi}(F) \ar[r]^-{i^{*}} &
\underline{\mathbf{G}}^{\chi}(U)}$$

\item Pour tout champ alg\'ebrique $F$, il existe un morphisme $\underline{\mathbf{K}}^{\chi}(F) \longrightarrow 
\underline{\mathbf{G}}^{\chi}(F)$, fonctoriel pour les images r\'eciproques. 

Si $F$ est r\'egulier, alors le morphisme ci-dessus est un isomorphisme, et est compatible avec les produits. 

\item Pour tout torseur affine sur un champ alg\'ebrique $p : V \longrightarrow F$, le morphisme
$$p^{*} : \underline{\mathbf{G}}^{\chi}(F) \longrightarrow \underline{\mathbf{G}}^{\chi}(V)$$
est un isomorphisme.
\end{itemize}

\textit{Remarque:} Comme le morphisme naturel $\pi : C_{F}^{t} \longrightarrow F$ poss\`ede une section
naturelle (correspondant au sous-groupe trivial), il existe des d\'ecompositions canoniques
$$\underline{\mathbf{K}}_{*}^{\chi}(F)\simeq \underline{\mathbf{K}}_{*}(F)\oplus \underline{\mathbf{K}}_{*}^{\chi\neq 1}(F)$$
$$\underline{\mathbf{G}}_{*}^{\chi}(F)\simeq \underline{\mathbf{G}}_{*}(F)\oplus \underline{\mathbf{G}}_{*}^{\chi\neq 1}(F)$$

Supposons que $F$ soit un champ alg\'ebrique r\'egulier. Alors on peut montrer que le morphisme naturel
$$\mathbf{K}(F) \longrightarrow \mathbf{G}(F)$$
est un isomorphisme. En effet, comme $F$ est r\'egulier, 
tout complexe de faisceaux quasi-coh\'erents sur $F$, \`a cohomologie coh\'erente et born\'ee, est un complexe parfait. Et, comme $F$ est 
noeth\'erien, on peut v\'erifier comme dans \cite[Cor. $3.13$]{th}, 
que $\mathbf{G}(F)$ est isomorphe au spectre de $K$-th\'eorie de ces complexes. \\

Le th\'eor\`eme principal est le suivant, dont la preuve sera donn\'ee dans le prochain paragraphe.

\begin{thm}\label{t1}
Soit $F$ un champ alg\'ebrique. 
\begin{enumerate}

\item Il existe un morphisme de spectres en anneaux
$$\phi_{F} : \mathbf{K}(F) \longrightarrow \underline{\mathbf{K}}^{\chi}(F),$$
qui est fonctoriel (dans la cat\'egorie homotopique) pour les images r\'eciproques. 

Si $F$ est r\'egulier, le morphisme $\phi_{F}$ ci-dessus induit un isomorphisme de spectres en anneaux
$$\phi_{F} : \mathbf{K}(F)_{\mathbb{Q}}\simeq \mathbf{G}(F)_{\mathbb{Q}} \simeq \underline{\mathbf{K}}^{\chi}(F)
\simeq \underline{\mathbf{G}}^{\chi}(F),$$

\item Si le morphisme $\pi : C_{F}^{t} \longrightarrow F$ est de $Tor$-dimension finie, alors il existe un moprhisme de spectres 
$$\phi_{F} : \mathbf{G}(F) \longrightarrow \underline{\mathbf{G}}^{\chi}(F),$$
qui est fonctoriel (dans la cat\'egorie homotopique) pour les images r\'eciproques. 

De m\^eme, si $F$ est une gerbe sur un espace alg\'ebrique, alors le morphisme ci-dessus $\phi_{F}$ induit un isomorphisme
de spectres
$$\phi_{F} : \mathbf{G}(F)_{\mathbb{Q}} \simeq \underline{\mathbf{G}}^{\chi}(F),$$

\item Il existe un isomorphisme de spectres 
$$\psi_{F} : \mathbf{G}(F)_{\mathbb{Q}} \longrightarrow \underline{\mathbf{G}}^{\chi}(F),$$
qui est fonctoriel (dans la cat\'egorie homotopique) pour les images directes de morphismes propres et repr\'esentables. 

\end{enumerate}
\end{thm}

Pour des raisons qui deviendrons claires lors de leur construction, les morphismes $\phi$ et $\psi$ seront appel\'es respectivement, morphisme
de diagonalisation et morphisme d'induction. Il serait cependant plus juste d'appeler morphisme d'induction le morphisme inverse de $\psi$. \\

Il est int\'eressant de d\'ecanuler ce th\'eor\`eme dans le cas o\`u $F=[X/H]$ est le champ quotient d'un sch\'ema $X$ par un groupe fini $H$.
Dans ce cas on a d\'ej\`a vu que 
$C_{F}^{t}\simeq \coprod[U^{c}/N_{c}]$, o\`u la somme est prise sur les classes de conjugaisons de sous-groupes cycliques de $H$. 
Supposons pour commencer que chaque $U^{c}$ contienne les racines $m(c)$-\`eme de l'unit\'e, o\`u $m(c)$ est l'ordre de $c$. Ainsi, 
la restriction de $\Lambda$ \`a $U^{c}$ est un faisceau constant isomorphe \`a $\mathbb{Q}(\zeta_{m(c)})$. Il est muni d'une
action du groupe $N_{c}$, ce qui en fait un faisceau localement constant sur $[U^{c}/N_{c}]$, qui est la restriction de $\Lambda$.
On a alors
$$\mathbb{H}^{-*}([U^{c}/N_{c}]_{et},\underline{\mathbf{K}}\wedge K(\mathbb{Q}(\zeta_{m(c)}),0))\simeq 
\left( \mathbf{K}_{*}(U^{c})\otimes \mathbb{Q}(\zeta_{m(c)}) \right)^{N_{c}},$$
et donc, d'apr\`es le point $(1)$ du th\'eor\`eme, on dispose d'un morphisme de $\mathbb{Q}$-alg\`ebres
$$\phi_{F} : \mathbf{K}_{*}(F) \longrightarrow \underline{\mathbf{K}}_{*}^{\chi}(F)\simeq
\bigoplus_{c \in c(H)}\left(\mathbf{K}_{*}(U^{c})\otimes \mathbb{Q}(\zeta_{m(c)})\right)^{N_{c}}.$$

Supposons que $F$ soit r\'egulier. Alors comme $c$ est un sch\'ema en groupes de type multiplicatif sur l'ouvert de $X$ o\`u son ordre 
est inversible, le sch\'ema $U^{c}$ est encore r\'egulier (\cite[$6.2$]{th2}). Le th\'eor\`eme affirme alors l'existence d'un isomorphisme de $\mathbb{Q}$-alg\`ebres
$$\phi_{F} : \mathbf{K}_{*}(F)\otimes \mathbb{Q} \simeq \bigoplus_{c \in c(H)}
\left(\mathbf{K}_{*}(U^{c})\otimes \mathbb{Q}(\zeta_{m(c)})\right)^{N_{c}}.$$
Notons que cette formule est essentiellement la formule d\'emontr\'ee par A. Vistoli dans \cite{v1}.

Quand au point $(3)$, il affirme l'existence d'un isomorphisme de $\mathbb{Q}$-espaces vectoriels
$$\psi_{F} : \mathbf{G}_{*}(F)\otimes \mathbb{Q} \simeq \bigoplus_{c \in c(H)}
\left(\mathbf{G}_{*}(U^{c})\otimes \mathbb{Q}(\zeta_{m(c)})\right)^{N_{c}}.$$
C'est aussi essentiellement la formule d\'emontr\'ee dans \cite{v1}. \\

Comme localement sur son espace de modules tout champ alg\'ebrique est de la forme $[X/H]$, on peut l\'egitimement penser que le th\'eor\`eme
\ref{t1} est obtenue par recollement des formules pr\'ec\'edentes pour la topologie \'etale. Notons que c'est ce proc\'ed\'e de recollement
qui fait apparaitre la cohomologie g\'en\'eralis\'ee \`a valeurs dans les pr\'efaisceaux de $K$-th\'eorie ou de $G$-th\'eorie. Je dois
avouer ne pas connaitre de d\'emonstration de \ref{t1} qui n'utilise pas cette th\'eorie homotopique (m\^eme pour le cas des champs 
alg\'ebriques lisses sur $Spec \mathbb{C}$ par exemple). \\

Terminons par l'\'enonc\'e de la formule de Lefschetz-Riemann-Roch, reliant les morphismes $\phi$ et $\psi$. Nous en donnerons une
esquisse de preuve dans la section suivante. Pour les d\'etails, nous renvoyons \`a \cite[Thm. $3.25$]{t2}. 

Pour cela, notons $\pi : C_{F}^{t} \longrightarrow F$ la projection naturelle. Comme c'est un morphismes non-ramifi\'e, il poss\`ede un
fibr\'e conormal, $\mathcal{N}^{\vee}$. Consid\'erons $\lambda_{-1}(\mathcal{N}^{\vee}):=\sum_{i}(-1)^{i}[\Lambda^{i}(\mathcal{N}^{\vee})] 
\in \mathbf{K}_{0}(C_{F}^{t})$. On pose alors (voir la d\'efinition de $d_{F} : \mathbf{K}_{0}(C_{F}^{t}) \longrightarrow 
\underline{\mathbf{K}}_{0}^{\chi}(F)$ dans le paragraphe suivant)
$$\alpha_{F}:=d_{F}(\lambda_{-1}(\mathcal{N}^{\vee})).$$
Une \'etude locale montre que $\alpha_{F}$ est partout de rang non-nul, et donc inversible dans $\underline{\mathbf{K}}_{0}^{\chi}(F)$
(\cite[$4.7$]{t1}).

\begin{thm}\label{t2}
Pour tout champ alg\'ebrique $F$, les morphismes
$$\psi_{F} : \mathbf{G}_{*}(F) \longrightarrow \underline{\mathbf{G}}_{*}^{\chi}(F)$$
$$\phi_{F} : \mathbf{K}_{*}(F) \longrightarrow \underline{\mathbf{K}}_{*}^{\chi}(F),$$
satisfaisont aux conditions suivantes.
\begin{enumerate}

\item Le morphisme $\phi_{F}$ est un morphisme d'anneaux, fonctoriels pour les images r\'eciproques. 

\item Soit $p : F \longrightarrow F'$ un morphisme propre de dimension cohomologique finie, et 
$x \in \mathbf{G}_{m}(F)$. Alors, on a
$$\psi_{F'}\circ p_{*}(x)=p_{*}\circ \psi_{F}(x)$$
dans un des trois cas suivants
\begin{itemize}
\item Le morphisme $p$ est repr\'esentable.
\item Le champ $F$ est r\'egulier.
\item $m=0$.
\end{itemize}

\item Si $F$ est lisse et quasi-projectif, alors pour tout $x \in \mathbf{G}_{*}(F)\simeq \mathbf{K}_{*}(F)$, on a
$$\psi(x)=\alpha_{F}^{-1}.\phi_{F}(x).$$

\item Pour tout morphisme \'etale $f : F \longrightarrow F'$, et $x \in \mathbf{G}_{m}(F')$ on a 
$$f^{*}\circ \psi_{F'}(x) = \psi_{F}\circ f^{*}(x)$$
dans l'un des cas suivant
\begin{itemize}
\item Le morphisme $f$ est repr\'esentable.
\item Le champ $F'$ est r\'egulier.
\item $m=0$.
\end{itemize}

\item Pour tout $x \in \mathbf{G}_{*}(F)$ et $y \in \mathbf{K}_{*}(F)$, on a
$$\psi_{F}(y.x)=\phi_{F}(y).\psi_{F}(x)$$

\end{enumerate}
\end{thm}

\end{section}

\begin{section}{D\'emonstration du th\'eor\`eme principal}

Pour d\'emontrer les th\'eor\`emes \ref{t1} et \ref{t2}, nous aurons besoin des r\'esultats de descente pour la $G$-th\'eorie des champs alg\'ebriques,
ainsi que de la notion de quasi-enveloppe de Chow et de leur existence. Ces
r\'esultats sont d\'emontr\'es dans \cite{t2}, et nous ne ferons que les rappeler.

\begin{subsection}{Rappels sur les r\'esultats de descente}

La proposition suivante permet de localiser sur l'espace de modules, et ainsi de ramener de nombreux \'enonc\'es au cas 
des champs quotients par des groupes finis. 

\begin{prop}\label{p1}
Soit $p : F \longrightarrow M$ la projection d'un champ alg\'ebrique $F$ sur son espace de modules, et consid\'erons le pr\'efaisceau en spectres 
$$\begin{array}{cccc}
p_{*}\mathbf{G}_{\mathbb{Q}} : & M_{et} & \longrightarrow & Sp \\
& U & \mapsto & \mathbf{G}(p^{-1}(U))_{\mathbb{Q}}
\end{array}$$
Alors le morphisme naturel
$$\mathbf{G}(F)_{\mathbb{Q}} \longrightarrow \mathbb{H}(M_{et},p_{*}\mathbf{G}_{\mathbb{Q}})$$
est un isomorphisme.

En particulier, le pr\'efaisceau en spectres $p_{*}\mathbf{G}_{\mathbb{Q}}$ est flasque sur $M_{et}$.
\end{prop}

\textit{Preuve:} Voir \cite[Thm. $2.4$]{t2}. $\Box$ \\

Remarquons que la proposition pr\'ec\'edente est fausse si l'on remplace $p_{*}$ par $\mathbb{R}p_{*}$. En effet, il n'est pas vrai que
le morphisme naturel $\mathbf{G}(F)_{\mathbb{Q}} \longrightarrow \mathbb{H}(F_{et},\underline{\mathbf{G}}_{\mathbb{Q}})$ est un
isomorphisme. Ceci est la principale diff\'erence avec le cas des sch\'emas (ou m\^eme des espaces alg\'ebriques), o\`u 
la $G$-th\'eorie rationnelle satisfait \`a la propri\'et\'e de descente \'etale. Ce ph\'enom\`ene explique aussi pourquoi la
cohomologie usuelle des champs alg\'ebriques n'est pas un invariant assez fin pour poss\'eder des formules
de Riemann-Roch (voir la remarque qui suit \cite[$4.3$]{t1}). \\

Soit $X \longrightarrow F$ un morphisme d'un espace alg\'ebrique $X$ vers un champ alg\'ebrique. Nous noterons $\mathcal{N}(X/F)$ le nerf
de ce morphisme. C'est un espace alg\'ebrique simplicial, o\`u le '$m$-\`eme \'etage' est donn\'e par
$$\mathcal{N}(X/F)_{m}:=\underbrace{X\times_{F}X \dots \times_{F}X}_{m+1 \; fois}$$
et les faces et d\'eg\'en\'erescences sont induits par les projections et les diagonales.

Supposons que $X \longrightarrow F$ soit propre, et donc que les morphismes de transitions $\mathcal{N}(X/F)_{m} \longrightarrow
\mathcal{N}(X/F)_{p}$ soient aussi des morphismes propres. On dispose alors 
d'un spectre simplicial, $[m] \mapsto \mathbf{G}(\mathcal{N}(X/F)_{m})$ (ici il s'agit d'un objet dans la cat\'egorie homotopique
des spectres simpliciaux et non d'un objet simplicial dans la cat\'egorie homotopique !). La colimite homotopique de ce spectre
simplicial sera not\'ee
$$\mathbf{G}(X/F):=Hocolim_{[m] \in \Delta^{o}}\mathbf{G}(\mathcal{N}(X/F)_{m}).$$

\begin{prop}\label{p2}
Soit $p : X \longrightarrow F$ un morphisme propre, avec $X$ un espace alg\'ebrique. Alors, le morphisme naturel
$$p_{*} : \mathbf{G}(X/F)_{\mathbb{Q}} \longrightarrow \underline{\mathbf{G}}(F)$$
est un isomorphisme.
\end{prop}

\textit{Preuve:} Voir \cite[Thm. $3.9$]{t1}. $\Box$ \\

\begin{cor}\label{c1}
Pour tout champ alg\'ebrique d'espace de modules $M$, il existe une isomorphisme covariant pour les
morphismes propres
$$\underline{\mathbf{G}}(F) \longrightarrow \mathbf{G}(M)_{\mathbb{Q}}.$$
\end{cor}

\textit{Preuve:} Voir \cite[$3.11$]{t1}. $\Box$ \\

Pour la d\'efinition suivante, rappelons qu'un point d'un champ $F$ est un point de son espace de modules $M$ (ceci n'est pas la d\'efinition
adopt\'ee dans \cite{lm}). De plus, pour $x \in M$ un point de $F$, on dispose de la gerbe r\'esiduelle 
en $x$, not\'ee $\widetilde{x}$. C'est une gerbe sur $Spec k(x)$, li\'ee par le groupe d'isotropie de $x$, $H_{x}$ (qui est bien 
d\'efini \`a conjugaison pr\`es). 

\begin{df}\label{d4}
Soit $p : F \longrightarrow F'$ un morphisme propre (non n\'ecessairement repr\'esentable) 
de champs alg\'ebriques. On dira que $p$ est une quasi-enveloppe de Chow
si pour tout point $x$ de $F'$, il existe un point $y$ dans $F$ avec $f(y)=x$, et tel que
le morphisme induit
$$p : \widetilde{y} \longrightarrow \widetilde{x}$$
admette une section apr\`es une extension finie du corps $k(x)$. 
\end{df}

Il est important de remarquer que l'existence de la section implique que le morphisme ci-dessus 
induit une surjection $p : H_{y} \longrightarrow H_{x}$. Ainsi, si $p : F \longrightarrow F'$ est une enveloppe, 
le morphisme induit $Ip : I_{F} \longrightarrow I_{F'}$ sur les champs des ramifications est surjectif (i.e. induit 
une surjection sur les espaces de modules). Ceci montre par exemple que le morphisme naturel $X \longrightarrow [X/H]$, pour
$H$ est un groupe fini, est une quasi-enveloppe de Chow si et seulement si $H$ est trivial. \\

Soit $F' \longrightarrow F$ un morphisme de champs alg\'ebriques. Alors le nerf de ce morphisme est un champ alg\'ebrique simplicial
augment\'e vers $F$, not\'e $\mathcal{N}(F'/F) \longrightarrow F$. Remarquons que ce n'est pas un objet simplicial de la cat\'egorie
homotopique des champs, mais un objet dans la cat\'egorie homotopique des champs simpliciaux. Si $F' \longrightarrow F$ est un 
morphisme propre (disons repr\'esentable), alors on peut appliquer le foncteur covariant $\mathbf{G}$, et ontenir un
spectre simplicial augment\'e, $[m] \mapsto \mathbf{G}(\mathcal{N}(F'/F)_{m})$. La colimite homotopique de ce spectre simplicial
sera not\'ee $\mathbf{G}(F'/F)$. 

\begin{prop}\label{p3}
Soit $p : F' \longrightarrow F$ un morphisme propre repr\'esentable, qui est une quasi-enveloppe de Chow. Alors, le morphisme naturel
$$p_{*} : \mathbf{G}(F'/F)_{\mathbb{Q}} \longrightarrow \mathbf{G}(F)_{\mathbb{Q}}$$
est un isomorphisme.
\end{prop}

\textit{Preuve:} Voir \cite[Prop. $3.13$]{t1}. $\Box$ \\

\end{subsection}

\begin{subsection}{Existence des enveloppes}

Pour la d\'efinition suivante nous noterons $\Delta_{+}$ la cat\'egorie des ensembles (\'eventuellement vides) 
finis totalement ordonn\'es. Les foncteurs $\Delta_{+}^{o} \longrightarrow C$ sont exactement les objets 
simpliciaux augment\'es dans $C$. \\

Fixons-nous un champ alg\'ebrique $F_{0}$, et consid\'erons la cat\'egorie des champs alg\'ebriques 
repr\'esentables sur $F_{0}$. Elle a pour objet les $1$-morphismes repr\'esentables 
$F \longrightarrow F_{0}$. Un morphisme entre $p : F \longrightarrow F_{0}$ et $q : F' \longrightarrow F_{0}$ est la
donn\'ee d'un $1$-morphisme repr\'esentable $f : F \longrightarrow F'$ et d'un $2$-morphisme $h : q\circ f \Rightarrow p$. Notons cette
cat\'egorie $Ch/F_{0}$. Par d\'efinition, un champ simplicial augment\'e sur $F_{0}$ est un objet simplicial de $Ch/F_{0}$. 

\begin{df}\label{d5}
Soit $F_{\bullet} \longrightarrow F$ un champ simplicial augment\'e. On dira que $F_{\bullet}$ est une enveloppe de $F$, si 
les propri\'et\'es suivantes sont satisfaites
\begin{enumerate}
\item Pour tout $[m] \in \Delta^{o}$, le champ $F_{m}$ est une r\'eunion disjointe de gerbes triviales sur des sch\'emas quasi-projectifs 
(i.e. est \'egal \`a une somme de $X_{m}\times_{S} BH_{m}$, avec 
$H_{m}$ un groupe fini et $X_{m}$ un sch\'ema quasi-projectif). 
\item Pour tout morphisme dans $\Delta^{o}_{+}$, $[m] \rightarrow [p]$, le morphisme $F_{m} \longrightarrow F_{p}$ est 
une somme de morphismes de la forme $f\times \rho : X_{m}\times_{S}BH_{m} \longrightarrow X_{p}\times BH_{p}$, avec $f : X_{m} \longrightarrow
X_{p}$ un morphisme propre, et $\rho$ un morphisme injectif de groupes (on parlera de morphismes cart\'esiens entre
gerbes triviales). 
\item Pour tout $[m] \in \Delta^{o}$, le morphisme naturel 
$$F_{m} \longrightarrow Cosq_{m-1}^{F}Sq_{m-1}(F_{\bullet})_{m}$$
est une quasi-enveloppe de Chow.
\end{enumerate}
\end{df}

Comme pr\'ec\'edemment, si $F_{\bullet} \longrightarrow F$ est une enveloppe, on dispose du spectre simplicial $[m] \mapsto \mathbf{G}(F_{m})$, 
dont la colimite homotopique sera not\'ee $\mathbf{G}(F_{\bullet})$. 

\begin{prop}\label{p4}
Soit $F$ un champ alg\'ebrique, et $p : F_{\bullet} \longrightarrow F$ une enveloppe. Alors le morphisme naturel
$$p_{*} : \mathbf{G}(F_{\bullet})_{\mathbb{Q}} \longrightarrow \mathbf{G}(F)_{\mathbb{Q}}$$
est un isomorphisme.
\end{prop}

\textit{Preuve:} Ceci provient de \ref{p3}, et d'un principe g\'en\'eral (voir \cite[Thm. $4.1$]{g} par exemple). $\Box$\\

Le principal r\'esultat d'existence est le suivant.

\begin{thm}\label{t3}
\begin{enumerate}
\item
Tout champ alg\'ebrique poss\`ede une enveloppe. 
\item
Si $F_{0}$ est un champ alg\'ebrique, et $F_{\bullet} \longrightarrow F_{0} \longleftarrow F_{\bullet}'$
deux enveloppes, il existe un diagramme commutatif
$$\xymatrix{
F_{\bullet}'' \ar[r] \ar[d] & F_{\bullet} \ar[d] \\
F_{\bullet}' \ar[r] & F_{0}}$$
o\`u $F_{\bullet}'' \longrightarrow F$ est une enveloppe. De plus, on peut choisir $F_{\bullet}''$ tel que pour tout $[m] \in \Delta^{o}$, 
les morphismes $F_{m}''\longrightarrow F_{m}$ et $F_{m}'' \longrightarrow F_{m}'$ soient cart\'esiens.
\item Si $F_{0} \longrightarrow F_{0}'$ est un morphisme repr\'esentable et propre de champs alg\'ebriques, il existe un diagramme
commutatif
$$\xymatrix{
F_{\bullet} \ar[r] \ar[d]^{p} & F_{\bullet}' \ar[d] \\
F_{0} \ar[r] & F_{0}'}$$
avec $F_{\bullet} \longrightarrow F_{0}$ et $F_{\bullet}' \longrightarrow F_{0}'$ des enveloppes. De plus, on peut choisir
$p$ de telle sorte que $F_{m} \longrightarrow F_{m}'$ soit cart\'esien pour tout $[m] \in\Delta^{o}$.
\end{enumerate}
\end{thm}

\textit{Preuve:} Comme la notion de quasi-enveloppe de Chow est stable par changements de bases et par compositions, il 
nous suffit de montrer que tout champ alg\'ebrique $F$ admet une quasi-enveloppe de Chow $F' \longrightarrow F$. En effet, 
la construction d'une enveloppe pour $F$ provient alors des m\'ethodes de construction de \cite[$V^{bis}.5.1.7$]{sga4}. 
Les assertions concernant l'aspect cart\'esien des morphismes provient du fait qu'un morphisme repr\'esentable
entre gerbes triviales $X\times_{S}BH \longrightarrow Y\times_{S}BK$, devient cart\'esien apr\`es un changement de base
fini et \'etale de $X$. Ainsi, dans la construction de l'enveloppe, on peut par induction rendre tous les morphismes de faces
cart\'esiens. De plus, comme les enveloppes sont des champs simpliciaux $\sigma$-scind\'es, les d\'eg\'en\'erescences
sont des facteurs directs, et donc clairement cart\'esiens. \\

Le champ $F$ \'etant de type fini sur $S$, il provient par changement de base d'un champ de type fini sur $Spec \mathbb{Z}$. Comme
les quasi-enveloppes de Chow sont stables par changements de bases, il s'en suit que l'on peut supposer que
$S=Spec \mathbb{Z}$. De plus, comme l'immersion ferm\'ee $F_{red} \hookrightarrow F$ est une
quasi-enveloppe de Chow, on peut aussi supposer que $F$ est r\'eduit. On peut aussi clairement supposer que $F$ est int\`egre.\\

En proc\'edant par r\'ecurrence noeth\'erienne, il nous suffit de contruire $F' \longrightarrow F$, propre et surjectif, qui est
g\'en\'eriquement une quasi-enveloppe de Chow. 

Comme $F$ est g\'en\'eriquement une gerbe sur un sch\'ema, le morphisme de normalisation $F' \longrightarrow F$
est g\'en\'eriquement une quasi-enveloppe. On peut donc supposer que $F$ est normal.

Soit $F \longrightarrow M$ l'espace de modules de $M$, et $X \longrightarrow M$ un morphisme g\'en\'eriquement
fini, propre et surjectif, avec $X$ un sch\'ema normal et quasi-projectif. Par \cite[Thm. $2.1$]{e}, 
on peut m\^eme supposer qu'il existe un diagramme commutatif
$$\xymatrix{ X \ar[r] \ar[rd] & F \ar[d] \\
& M}$$
Soit $F_{0}:=(F\times_{M}X)_{red} \longrightarrow F$, et $F'\longrightarrow F$ la normalisation de $F_{0}$. Le morphisme
$F' \longrightarrow F$ est clairement g\'en\'eriquement une quasi-enveloppe de Chow. Il nous reste donc \`a montrer que
$F'$ est une gerbe sur $X$. 

Notons $M'$ l'espace de modules de $F'$. Alors $M' \longrightarrow X$ est un morphisme fini et birationnel. Comme $M'$ et 
$X$ sont normaux, c'est un isomorphisme. Ainsi, l'espace de modules de $F'$ est $X$. De plus, par construction, le morphisme
naturel $F' \longrightarrow X$ poss\`ede une section. 

\begin{lem}\label{l1}
Soit $F$ un champ alg\'ebrique int\`egre et normal. Si la projection sur son espace de modules $F \longrightarrow M$
poss\`ede une section, alors $F$ est une gerbe sur $M$.
\end{lem}

\textit{Preuve:} Comme l'assertion est locale pour la topologie \'etale sur $M$, on peut supposer que $F=[X/H]$ est un champ
quotient d'un sch\'ema normal et int\`egre $X$ par l'action d'un groupe fini $H$. L'existence de la section
$s : M=X/H \longrightarrow F$ se traduit par un diagramme commutatif
$$\xymatrix{Y \ar[r]^{f} \ar[rd]_{q} & X \ar[d]^{p} \\
& X/H}$$
o\`u $q$ est un $H$-torseur, et $f$ est $H$-\'equivariant. Comme $q$ est \'etale, $f$ est non-ramifi\'e. Or $Y$ et $X$ sont
int\`egres normaux et de m\^eme dimension, et donc $f$ est \'etale. Ceci implique que $p$ est \'etale, et donc
que l'action se factorise par une action libre de $H/H_{0}$, o\`u $H_{0}$ est le noyau de l'action de $H$ sur $X$. $\Box$ \\

Ce lemme implique que $F'$ est une gerbe sur $X$. De plus, comme cette gerbe est neutre, elle est de la forme $BH$, o\`u 
$H \longrightarrow X$ est un sch\'ema en groupes fini et \'etale sur $X$. Ainsi, en effectuant un changement de base
fini et \'etale de $X$, on peut supposer que $F'\simeq X\times BH$. Ce qui ach\`eve la d\'emonstration. $\Box$\\

\textit{Remarque:} La preuve pr\'ec\'edente montre aussi que pour tout champ alg\'ebrique $F$, il existe un
morphisme repr\'esentable et fini $F_{0} \longrightarrow F$, qui est une quasi-enveloppe de Chow, et avec $F_{0}=\coprod X\times_{S}BH$. 

\end{subsection}

\begin{subsection}{D\'emonstration du th\'eor\`eme}

\textit{Construction de $\phi_{F}$:} \\

Rappelons que l'on dispose du faisceau des caract\`eres $\mathcal{X}$ sur $(C_{F}^{t})_{et}$, et qu'il s'agit d'un faisceau 
fini et localement constant pour la topologie \'etale. On dispose aussi du pr\'efaisceau en cat\'egories ab\'eliennes
$$\begin{array}{cccc}
QCoh : & (C_{F}^{t})_{et} & \longrightarrow & CatAb \\
& U & \mapsto & QCoh(U)
\end{array}$$
Notons $\prod_{\mathcal{X}}QCoh$ le pr\'efaisceau en cat\'egories exactes d\'efini par
$$\begin{array}{cccc}
\prod_{\mathcal{X}}QCoh : & (C_{F}^{t})_{et} & \longrightarrow & CatAb \\
& U & \mapsto & (\prod_{\mathcal{X}}QCoh)(U):=\prod_{\mathcal{X}(U)}QCoh(U)
\end{array}$$
Remarquons que ce pr\'efaisceau en cat\'egories n'est pas un champ. Le champ qui lui est associ\'e sera not\'e
$\mathcal{X}\otimes QCoh$. La cat\'egories des sections cart\'esiennes de ce champ sur $(C_{F}^{t})_{et}$ sera not\'ee
$\mathcal{X}\otimes QCoh(C_{F}^{t})$. C'est naturellement une cat\'egorie ab\'elienne.\\

Commen\c{c}ons par construire un foncteur exact $d_{F} : QCoh(C_{F}^{t}) \longrightarrow \mathcal{X}\otimes QCoh(C_{F}^{t})$, 
construit par diagonalisation de l'action du champ en groupes universel $\mathcal{C}^{t}_{F} \longrightarrow C_{F}^{t}$.

Pour cela, remarquons que tout faisceau quasi-coh\'erent sur $C_{F}^{t}$ arrive avec une
action naturelle de $\mathcal{C}_{F}^{t}$. 
Ceci d\'efinit un foncteur exact 
$$QCoh(C_{F}^{t}) \longrightarrow QCoh(C_{F}^{t},\mathcal{C}^{t})$$
de la cat\'egorie des faisceaux quasi-coh\'erents sur $C_{F}^{t}$ vers celle des faisceaux quasi-coh\'erents sur $C_{F}^{t}$ munis d'une action
de $\mathcal{C}_{F}^{t}$. On utilise alors le lemme suivant.

\begin{lem}\label{l2}
Soit $F$ un champ alg\'ebrique, et $C \longrightarrow F$ un champ en groupes fini de type multiplicatif. Notons $\mathcal{X}$ son
faisceau des caract\`eres sur $F_{et}$. Alors il existe une \'equivalence naturelle de cat\'egories ab\'eliennes
$$QCoh(F,C) \simeq \mathcal{X}\otimes QCoh(F)$$
o\`u $QCoh(F,C)$ est la cat\'egorie des faisceaux quasi-coh\'erents $C$-\'equivariants sur $F$, et $\mathcal{X}\otimes QCoh(F)$ la
cat\'egorie des sections globales cart\'esiennes du champ associ\'e au pr\'efaisceau en cat\'egories sur $F$,
$U \mapsto \prod_{\mathcal{X}(U)}QCoh(U)$.
\end{lem}

\textit{Preuve:} Soit $V$ un faisceau quasi-coh\'erent $C$-\'equivariant sur $F$, et $U \longrightarrow F$ un morphisme \'etale. 
La restriction de $V$ \`a $U$ poss\`ede une action du sch\'ema en groupes fini $C_{U}$. Choisissons $U$ tel que
$C_{U}$ soit constant sur $U$. Alors le faisceau se d\'ecompose en une somme directe
$$V_{U}\simeq \bigoplus_{\chi \in \mathcal{X}(U)}V_{U}^{(\chi)}$$
o\`u $V_{U}^{(\chi)}$ est le sous-faisceau o\`u $C_{U}$ op\`ere par le caract\`ere $\chi$. On consid\`ere alors
$\prod_{\chi \in \mathcal{X}(U)}V_{U}^{(\chi)} \in \prod_{\mathcal{X}(U)}QCoh(U)\simeq \mathcal{X}\otimes QCoh(U)$.
On peut v\'erifier que les cocycles de $V$ sur $F_{et}$ induisent des cocycles pour les 
$\prod_{\chi \in \mathcal{X}(U)}V_{U}^{(\chi)}$ lorsque $U$ varie dans $F_{et}$. Ils se recollent donc un objet bien d\'efini
dans $\mathcal{X}\otimes QCoh(F)$. On laisse le soin au lecteur de v\'erifier que ceci d\'efinit bien un foncteur exact
$$QCoh(F,C) \longrightarrow \mathcal{X}\otimes QCoh(F).$$
Pour montrer que c'est une \'equivalence, on utilise que $U \mapsto QCoh(U,C_{U})$ et $U \mapsto \mathcal{X}\otimes QCoh(U)$
sont tous deux des champs sur $F_{et}$. On peut donc remplacer $F$ par n'importe quel $U \longrightarrow F$
qui est \'etale et surjectif. On peut donc supposer que $F$ est un sch\'ema, et que $C$ est diagonalisable. Le
r\'esultat est alors imm\'ediat. $\Box$\\

On vient donc de construire un foncteur exact $d_{F} : QCoh(C_{F}^{t}) \longrightarrow \mathcal{X}\otimes QCoh(C_{F}^{t})$. Ce
foncteur induit donc un morphisme sur les spectres de $K$-th\'eorie des complexes parfaits de faisceaux quasi-coh\'erents
$$d_{F} : \mathbf{K}(C_{F}^{t}) \longrightarrow \mathbf{K}(\mathcal{X}\otimes C_{parf}QCoh(C_{F}^{t})).$$
On utilise alors le fait suivant (voir aussi \cite[Prop. $1.6$]{t2}). Soit $C$ un site, muni d'un pr\'efaisceau en cat\'egories compliciale
de Waldhausen $\mathcal{E}$ (\cite{th}),
de section globales cart\'esiennes $\mathcal{E}(C)$. Notons $\underline{\mathbf{K}}$ le pr\'efaisceau en spectres
$U \mapsto \mathbf{K}(\mathcal{E}(U))$. Alors il existe un morphisme naturel
$$\mathbf{K}(\mathcal{E}(C)) \longrightarrow \Gamma(C,\underline{\mathbf{K}}) \longrightarrow 
\mathbb{H}(C,\underline{\mathbf{K}}).$$
Pour nous, ceci implique qu'il existe un morphisme naturel
$$\mathbf{K}(\mathcal{X}\otimes C_{parf}QCoh(C_{F}^{t})) \longrightarrow \mathbb{H}((C_{F}^{t})_{et},\mathcal{X}\otimes \underline{\mathbf{K}}).$$
Or, $(\mathcal{X}\otimes \underline{\mathbf{K}})_{\mathbb{Q}}\simeq \underline{\mathbf{K}}\wedge K(\mathbb{Q}[\mathcal{X}],0)$. On
peut donc composer avec la projection naturelle $\mathbb{Q}[\mathcal{X}] \longrightarrow \Lambda$ pour obtenir un morphisme
$$\mathbf{K}(\mathcal{X}\otimes C_{parf}QCoh(C_{F}^{t}))_{\mathbb{Q}} \longrightarrow 
\mathbb{H}((C_{F}^{t})_{et},\mathcal{X}\otimes \underline{\mathbf{K}}_{\mathbb{Q}})
\longrightarrow \mathbb{H}((C_{F}^{t})_{et},\underline{\mathbf{K}}_{\Lambda}).$$
En conclusion, on trouve un morphisme de spectres
$$d_{F} : \mathbf{K}(C_{F}^{t}) \longrightarrow \underline{\mathbf{K}}^{\chi}(F).$$
Il est facile de v\'erifier que ce morphisme est un morphisme de spectres en anneaux. \\

On d\'efinit $\phi_{F}$ par la composition
$$\phi_{F} : \xymatrix{
\mathbf{K}(F) \ar[r]^{\pi^{*}} & \mathbf{K}(C_{F}^{t}) \ar[r]^-{d_{F}} & \underline{\mathbf{K}}^{\chi}(F)}$$
o\`u $\pi : C_{F}^{t} \longrightarrow F$ est la projection naturelle. \\

Pour montrer que le morphisme $\phi_{F}$ est un isomorphisme lorsque $F$ est r\'egulier, on peut utiliser \ref{p1} (car
$\mathbf{G}(F)\simeq \mathbf{K}(F)$), et donc localiser (pour la topologie
\'etale) le probl\`eme sur l'espace de modules $M$ de $F$. On peut donc supposer que $F=[X/H]$ est un champ quotient
d'un sch\'ema par un groupe fini. En effectuant un changement de base par un morphisme \'etale $M' \longrightarrow M=X/H$, on peut
aussi supposer que $X$ contient les racines $m$-\`eme de l'unit\'e, o\`u $m$ est l'ordre de $H$. On se ram\`ene donc au cas o\`u
$$\underline{\mathbf{K}}_{*}^{\chi}(F)\simeq \bigoplus_{c \in c(H)}\left(\mathbf{K}_{*}(U^{c})\otimes \mathbb{Q}(\zeta_{m(c)})\right)^{N_{c}}.$$
Le morphisme $\mathbf{K}_{*}(F) \longrightarrow \left(\mathbf{K}_{*}(U^{c})\otimes \mathbb{Q}(\zeta_{m(c)})\right)^{N_{c}}$ est alors 
construit de la fa\c{c}on suivante. On dipose alors d'un morphisme de champ
$$j_{c} : [U^{c}/N_{c}] \longrightarrow [X/H].$$
Ce morphisme induit un morphisme en $K$-th\'eorie
$$j_{c}^{*} : \mathbf{K}_{*}(F) \longrightarrow \mathbf{K}_{*}([U^{c}/N_{c}]).$$
Par restriction de $N_{c}$ \`a $c$, il existe aussi un morphisme naturel $\mathbf{K}_{*}([U^{c}/N_{c}]) \longrightarrow
\mathbf{K}_{*}([U^{c}/c])^{N_{c}}$. Or, comme $c$ op\`ere trivialement sur $U^{c}$ et que $c$ est diagonalisable sur $U^{c}$, on 
obtient $\mathbf{K}_{*}([U^{c}/c])_\mathbb{Q}\simeq \mathbf{K}_{*}(U^{c})\otimes \mathbb{Q}[c^{*}]$, o\`u $c^{*}=Hom_{gp}(c,\mathbb{C}^{*})$.
En utilisant la projection naturel $\mathbb{Q}[c^{*}] \longrightarrow \mathbb{Q}(\zeta_{m(c)})$, on en d\'eduit le morphisme
cherch\'e
$$\mathbf{K}_{*}(F) \longrightarrow \left(\mathbf{K}_{*}(U^{c})\otimes \mathbb{Q}(\zeta_{m(c)})\right)^{N_{c}}.$$
Pour montrer que la somme de ces morphismes
$$\mathbf{K}(F)_{\mathbb{Q}} \longrightarrow \bigoplus_{c \in c(H)}\left(\mathbf{K}_{*}(U^{c})\otimes \mathbb{Q}(\zeta_{m(c)})\right)^{N_{c}}$$
est un isomorphisme  on proc\`ede exactement comme dans \cite{v1}, avec de modestes changements pour tenir compte du fait que l'on est en caract\'eristiques mixtes. \\

\textit{D\'emonstration de $(2)$:} \\

Nous ne construirons pas $\phi_{F}$. En effet, il suffit de r\'ep\'eter la construction de $d_{F}$ mais en rempla\c{c}ant les
faisceaux quasi-coh\'erents par des faisceaux coh\'erents. On obtient de cette fa\c{c}on un foncteur exact
$$d_{F} : Coh(C_{F}^{t}) \longrightarrow \mathcal{X}\otimes Coh(C_{F}^{t}).$$
Par les m\^emes arguments que pr\'ec\'edemment, ce foncteur induit un morphisme de spectres
$$\mathbf{G}(C_{F}^{t}) \longrightarrow \underline{\mathbf{G}}^{\chi}(F).$$
En composant avec $\pi^{*} : \mathbf{G}(F) \longrightarrow \mathbf{G}(C_{F}^{t})$, qui existe car $\pi$ est 
de $Tor$-dimension finie, on obtient
$$\phi_{F} : \mathbf{G}(F) \longrightarrow \underline{\mathbf{G}}^{\chi}(F).$$
 
Lorsque $F$ est une gerbe, on peut de nouveau invoquer \ref{p1}, et donc supposer que 
$F=[X/H]$, o\`u $H$ op\`ere trivialement sur $X$. De plus, $\pi$ est alors \'etale, et le morphisme
$\phi_{F}$ compatible avec les suites exactes longues de localisation. Par d\'evissage on se ram\`ene
donc au cas o\`u $X$ est le spectre d'un corps contenant suffisement de racines de l'unit\'e. La formule dans ce cas est 
alors bien connue. \\

\textit{Construction de $\psi_{F}$:} \\

Le lemme fondamental est le suivant.

\begin{lem}\label{l3}
Notons $Gbt/F$, la sous-cat\'egorie des champs propres et repr\'esentables sur $F$, qui sont des gerbes triviales, et des $1$-morphismes
propres et cart\'esiens.
Alors, il existe un isomorphisme $\phi$ (dans la cat\'egorie homotopique des pr\'efaisceaux en spectres
sur $(Gbt/F)^{o}$) entre les deux pr\'efaisceaux en spectres suivants 
$$\begin{array}{ccc}
(Gbt/F)^{o} & \longrightarrow & Sp \\
F' & \mapsto & \mathbf{G}(F')_{\mathbb{Q}} \\
F' & \mapsto & \underline{\mathbf{G}}^{\chi}(F')
\end{array}$$
De plus, pour tout $F'\in Gbt/F$, le morphisme $\phi : \mathbf{G}(F') \longrightarrow \underline{\mathbf{G}}^{\chi}(F')$ coincide avec
le morphisme $\phi_{F'}$ de $(2)$ \ref{t1}.
\end{lem}

\textit{Preuve:} Comme les gerbes sont repr\'esentables sur $F$, leurs groupes d'isotropie sont des sous-groupes de ceux de $F$. Ainsi, il
existe un entier $n$ qui est un multiple de tous les ordres des groupes $H$ liant les gerbes de $Gbt/F$.  Soit $m$ le plus grand entier
divisant $n$ et premier avec la caract\'eristique de $K(S)$ (le corps des fonctions de $S$). Soit $S' \longrightarrow S$ le sch\'ema
en groupes des racines $m$-\`eme de l'unit\'e.  Alors, en effectuant une descente galoisienne de $S'$ \`a $S$, on peut supposer que 
$\mu_{m}$ est un faisceau constant sur $S$. \\

Soit $X\times_{S}H \in Gbt/F$, et $c$ un sous-groupe cyclique de $H$. Notons $X_{c}$ le sous-sch\'ema o\`u l'ordre de $c$ est inversible. 
Comme $S$ contient suffisement de racines de l'unit\'e, $X_{c}$ aussi, et $c$ est diagonalisable sur $X_{c}$. Ainsi, son faisceau
des caract\`eres est constant \'egal \`a $c^{*}\simeq \mathbf{Z}/m_{c}$. Notons $\Lambda_{c}$ le quotient $\mathbb{Q}[c^{*}] \longrightarrow
\Lambda$ correspondant \`a $\mathbb{Q}[\mathbb{Z}/m_{c}] \longrightarrow \mathbb{Q}(\zeta_{m_{c}})$. 

Soit $N_{c}$ le normalisateur de $c$ dans $H$. Ce groupe op\`ere par conjugaison sur $\Lambda_{c}$. Il est alors facile de v\'erifier qu'il
existe un \'equivalence faible naturelle
$$\alpha : \prod_{c \in c(H)}(\mathbf{G}(X_{c})\otimes \Lambda_{c})^{N_{c}} \longrightarrow \underline{\mathbf{G}}^{\chi}(X\times_{S}BH)$$
(ici ($\mathbf{G}(X_{c})\otimes \Lambda_{c})^{N_{c}}:=Holim_{N_{c}}(\mathbf{G}(X_{c})\otimes \Lambda_{c}$)).

Pour $f\times \rho : X\times_{S}BH \longrightarrow Y\times_{S}BK$ un morphisme de $Gbt/F$, on d\'efinit une image directe
$$\begin{array}{cccc}
(f\times \rho)_{*} : & (\mathbf{G}(X_{c})\otimes \Lambda_{c})^{N_{c}} & \longrightarrow & 
(\mathbf{G}(Y_{\rho(c)})\otimes \Lambda_{\rho(c)})^{N_{\rho(c)}} \\
& \mathcal{F}\otimes x & \mapsto & f_{*}(\mathcal{F})^{[K:H]}\otimes \rho(x)
\end{array}$$
En d'autres termes, on a $(f\times\rho)_{*}:=([K:H]\times f_{*})\otimes \rho$. Remarquons que $\rho : \Lambda_{c} \longrightarrow \Lambda_{\rho(c)}$
est un isomorphisme.

Ceci permet de construire une image directe
$$(f\times\rho)_{*} :  \prod_{c \in c(H)}(\mathbf{G}(X_{c})\otimes \Lambda_{c})^{N_{c}} \longrightarrow
\prod_{c \in c(K)}(\mathbf{G}(Y_{c})\otimes \Lambda_{c})^{N_{c}}$$
(par d\'efinition, les composantes de $(f\times \rho)_{*}$ sont nulles pour $c \in c(K)$ qui n'est pas l'image de $c \in c(H)$). 

On peut alors v\'erifier que l'\'equivalence $\alpha$ est fonctorielle pour les images directes de morphismes de $Gbt/F$. On identifie ainsi
$X\times_{S}BH \mapsto \underline{\mathbf{G}}^{\chi}(X\times_{S}BH)$ \`a $X\times_{S}BH \mapsto 
\prod_{c \in c(H)}(\mathbf{G}(X_{c})\otimes \Lambda_{c})^{N_{c}}$. \\

On construit alors $\mathbf{G}(X\times_{S}BH) \longrightarrow \prod_{c \in c(H)}(\mathbf{G}(X_{c})\otimes \Lambda_{c})^{N_{c}}$
exactement par le m\^eme proc\'ed\'e de construction que celui utilis\'e pour d\'efinir $\phi$ dans le point $(1)$ de \ref{t1}. C'est alors
un exercice de v\'erifier que ce morphisme est covariant pour les images directes de morphismes de $Gbt/F$ (il s'agit d'utiliser que
le caract\`ere d'une repr\'esentation induite est \'egal au caract\`ere induit). $\Box$\\

Si $F_{\bullet} \longrightarrow F$ est une enveloppe, la transformation $\phi$ permet de d\'efinir un morphisme de spectres simpliciaux, entre
$[m] \mapsto \mathbf{G}(F_{m})$ et $[m] \mapsto \underline{\mathbf{G}}^{\chi}(F_{m})$. En passant \`a la colimite homotopique, on obtient
un morphisme
$$\mathbf{G}(F_{\bullet}) \longrightarrow \underline{\mathbf{G}}^{\chi}(F_{\bullet}).$$
De plus, le deuxi\`eme membre poss\`ede naturellement un morphisme vers $\underline{\mathbf{G}}^{\chi}(F)$. On obtient ainsi
$$\mathbf{G}(F_{\bullet})_{\mathbb{Q}} \longrightarrow \underline{\mathbf{G}}^{\chi}(F),$$
et en utilisant \ref{p4}
$$\psi_{F} : \mathbf{G}(F)_{\mathbb{Q}} \longrightarrow \underline{\mathbf{G}}^{\chi}(F).$$
On peut v\'erifier que ce morphisme est ind\'ependant de l'enveloppe choisie, car deux telles enveloppes sont domin\'ees par une troisi\`eme.
De plus, par construction et par \ref{l3}, ce morphisme est clairement covariant pour les morphismes propres et repr\'esentables. \\

Pour montrer que $\psi_{F}$ est un isomorphisme, on peut proc\'eder par r\'ecurrence noeth\'erienne et se ramener, 
\`a l'aide de la suite exacte longue de localisation, au cas o\`u 
$F$ est une gerbe. De plus, par une localisation sur l'espace de modules on peut m\^eme supposer que $F$ est une gerbe triviale. On choisissant
l'enveloppe triviale pour $F$, on en d\'eduit que dans ce cas $\psi_{F}=\phi_{F}$, et donc par le point $(1)$, $\psi_{F}$ est un
isomorphisme. \\

\end{subsection}

\begin{subsection}{Formule de Lefschetz-Riemann-Roch}

Dans ce paragraphe nous donnerons une esquisse de preuve de \ref{t2}. \\

Commen\c{c}ons par remarquer que le th\'eor\`eme \ref{t2} est vrai pour les champs de la forme $[X/H]$, avec $H$ un groupe
fini op\'erant sur un sch\'ema quasi-projectif $X$. Ceci se v\'erifie directement \`a la main (c'est essentiellement la
formule de Lefschetz-Riemann-Roch pour l'action d'un groupe fini). \\

$(1)$ Il n'y a rien \`a d\'emontrer. \\

$(2)$ Soit $f : F \longrightarrow F'$ propre, et $x \in \mathbf{G}_{m}(F)$. Le cas o\`u $f$ est repr\'esentable est d\'ej\`a connu
par \ref{t1}. 

Supposons que $F$ soit r\'egulier. Notons $p : F_{0} \longrightarrow F$ une quasi-enveloppe de Chow, avec $p$ repr\'esentable fini,
et $F_{0}$ une gerbe triviale. Comme $F$ est r\'egulier, on dispose de la formule de projection 
$$p_{*}p^{*}(x)=p_{*}(1).x$$
On montre tout d'abord que $p_{*}(1)$ est inversible. 
Comme $F$ est r\'egulier, on dispose d'un isomorphisme d'anneaux (\ref{t1})
$$\phi_{F} : \mathbf{K}_{0}(F)_{\mathbb{Q}} \longrightarrow \underline{\mathbf{K}}_{0}^{\chi}(F).$$
Ainsi, le fait que $p_{*}(1)$ soit inversible dans $\mathbf{K}_{0}(F)_{\mathbb{Q}}$ est en fait une assertion locale
pour la topologie \'etale de l'espace de modules $M$. On peut
donc supposer que $F=[X/H]$. On utilise alors que $p_{*}\circ \psi_{F_{0}}=\alpha_{F}^{-1}.\phi_{F}\circ p_{*}$. Ainsi, il suffit de 
montrer que $p_{*}(1) \in \underline{\mathbf{K}}^{\chi}_{0}(F)$ est inversible, ce qui provient directement du fait que
le morphisme $C_{F_{0}}^{t} \longrightarrow C_{F}^{t}$ est fini et surjectif (car $p$ est une quasi-enveloppe de Chow), 
et donc que $p_{*}(1)$ est partout de rang non-nul sur $C_{F}^{t}$.

On vient donc de voir que $p_{*}(1)$ est inversible, et donc ceci entraine que $p_{*} : \mathbf{G}_{*}(F_{0})
\longrightarrow \mathbf{K}_{*}(F)_{\mathbb{Q}}$ est surjectif. Ceci permet de remplacer $F$ par $F_{0}$, et m\^eme de
supposer que $F$ et $F'$ sont des gerbes triviales (car on peut prendre $F_{0}$ qui domine une autre quasi-enveloppe de
Chow $F_{0}' \rightarrow F'$). Ce cas est alors un cas particulier du cas o\`u $F$ et $F'$ sont des quotients
par des groupes finis. 

Enfin, supposons $F$ quelconque, mais $x \in \mathbf{G}_{0}(F)$. Un raisonnement par r\'ecurrence sur la codimension du support
montre aussi que $p_{*} : \mathbf{G}_{0}(F_{0})  \longrightarrow \mathbf{G}_{0}(F)_{\mathbb{Q}}$ est surjectif. On peut alors refaire
l'argumentation pr\'ec\'edente. \\

$(3)$ Supposons que $F$ soit lisse (connexe) et quasi-projectif. Alors on peut \'ecrire $u : F \longrightarrow F_{0}$, o\`u $F_{0}$
est le champ quotient d'un sch\'ema quasi-projectif $X$ par $Gl_{m}$, et $u$ est \'etale et
fini (\cite[Cor. $2.3$]{e}). Soit $F' \longrightarrow F_{0}$ une quasi-enveloppe de Chow, avec $p$ repr\'esentable fini, et 
$F'$ une gerbe triviale. Comme $F_{0}$ est un quotient par $Gl_{m}$ d'un sch\'ema quasi-projectif, $p$ est un morphisme
fortement projectif (i.e. se factorise en une immersion ferm\'e suivie de la projection d'un fibr\'e projectif associ\'e \`a un 
fibr\'e vectoriel sur $F_{0}$). On consid\`ere $p : F'\times_{F_{0}}F \longrightarrow F$, qui est encore une quasi-enveloppe et
fortement projectif. De plus, $F'':=F'\times_{F_{0}}F$ \'etant une gerbe, on peut effectuer un changement de base fini de son espace 
de modules, et supposer que $F''$ est une gerbe triviale. Ainsi, on a trouv\'e $p : F'' \longrightarrow F$, qui est 
une composition de morphismes fortement projectifs, ainsi qu'une quasi-enveloppe de Chow, et avec $F''$ une gerbe triviale. 
Par les techniques standard (voir par exemple \cite[Lem. $4.2$]{t1} \'etape $(a)$) on peut montrer que
$$p_{*}\circ \phi_{F''}=\alpha_{F}^{-1}.\phi_{F}\circ p_{*}.$$
Comme de plus, comme $p_{*} : \mathbf{G}_{*}(F'') \longrightarrow \mathbf{G}_{*}(F)_{\mathbb{Q}}$ 
est surjectif, et que $p_{*}\circ \phi_{F''}=p_{*}\circ \psi_{F''}=\psi_{F}\circ p_{*}$, ceci montre que
$$\psi_{F}=\alpha_{F}^{-1}.\phi_{F}.$$

$(4)$ utilise les m\^emes m\'ethodes de r\'eduction au cas des gerbes triviales. \\

$(5)$ Soit $p : F_{\bullet} \longrightarrow F$ une enveloppe, et $x \in \mathbf{G}_{*}(F)_{\mathbb{Q}}$ et $y\in \mathbf{K}_{*}(F)_{\mathbb{Q}}$. 
Par \ref{p4}, on \'ecrit $x=p_{*}(z)$. Alors, par la formule de projection, on a $x.y=p_{*}(z.p^{*}(y))$. Notons
$\phi : \mathbf{G}_{*}(F_{\bullet})_{\mathbb{Q}} \longrightarrow \underline{\mathbf{G}}_{*}^{\chi}(C_{F_{\bullet}}^{t})$ le morphisme
construit dans la preuve du point $(3)$ de \ref{t1}. Alors, comme $\phi$ est compatible avec les produits et les images r\'eciproques, et par
la formule de projection, on a 
$$\phi(z.p^{*}(y))=\phi(z).p^{*}\phi(y).$$
Or, par d\'efinition, on $\psi_{F}(x.y)=\psi_{F}(p_{*}(x.p^{*}(y)))=p_{*}(\phi(z).p^{*}\phi_{F}(y))=p_{*}(\psi(x)).\phi_{F}(y)=\psi(x).\phi(y)$.
$\Box$ \\

\end{subsection}

\end{section}

\begin{section}{Deux applications}

\begin{subsection}{$K$-th\'eorie \'equivariante}

La premi\`ere application de \ref{t1} que nous proposons concerne les questions $(2.4)$, $(3.4)$ et $(3.6)$ 
pos\'ees dans \cite{v2}, ainsi que leurs g\'en\'eralisations au cas des champs de Deligne-Mumford. 

Signalons que la conjecture $(2.4)$ a d\'ej\`a \'et\'e d\'emontr\'ee en utilisant les groupes de Chow \'equivariants (\cite{eg}). De plus, 
dans un travail r\'ecent (\cite{vv}), A. Vistoli et G. Vezzosi donnent une description g\'en\'erale de la $K$-th\'eorie \'equivariante 
d'une action avec stabilisateurs finis (\'eventuellement non-r\'eduits), ce qui leur permet de r\'esoudre aussi les conjectures $(3.4)$ et $(3.6)$. 
Cette formule est donc plus g\'en\'erale que les formules que
nous donnons dans ce paragraphe, et les techniques de d\'emonstration tr\`es diff\'erentes de celle que nous utilisons.
Dans le cas o\`u les stabilisateurs sont r\'eduits, cette formule tensoris\'e par $\mathbb{Q}$ 
est essentiellement \'equivalente \`a la formule \ref{c6'}. D'un autre cot\'e, les corollaires \ref{c2} et \ref{c5}
restent vrais pour des champs de Deligne-Mumford qui ne sont pas n\'ecessairement des quotients par des groupes alg\'ebriques, et 
pour lesquels \cite{vv} ne peut plus s'appliquer. \\

Fixons nous $G$ un sch\'ema en groupes plat et de type fini sur $S$, op\'erant sur un sch\'ema $X$ (de type fini sur $S$). 
On suppose que l'action et propre, et que les stabilisateurs de cette action sont des sch\'emas en groupes finis et r\'eduits. 
Dans ce cas, le champ quotient $F=[X/G]$ est un champ alg\'ebrique de Deligne-Mumford s\'epar\'e sur $S$.

L'espace de modules de $F$ est donc $M=X/G$, et $\mathbf{G}_{*}(F)$ est alors la $G$-th\'eorie \'equivariante 
de $X$ (i.e. la $K$-th\'eorie de la cat\'egorie ab\'eliennne des $G$-faisceaux coh\'erents sur $X$).
Rappelons que si $X$ est un sch\'ema r\'egulier, alors le morphisme
$$\mathbf{K}_{*}(F) \longrightarrow \mathbf{G}_{*}(F)$$
est un isomorphisme.

Nous noterons $p : F \longrightarrow M$ la projection sur l'espace de modules. \\

\begin{cor}\label{c2}
Il existe un d\'ecomposition 
$$\mathbf{G}_{*}([X/G])_{\mathbb{Q}}\simeq 
\underline{\mathbf{G}}_{*}([X/G])_{\mathbb{Q}}\oplus \underline{\mathbf{G}}_{*}^{\chi\neq 1}([X/G])
\simeq \mathbf{G}_{*}(X/G)_{\mathbf{Q}}\oplus \underline{\mathbf{G}}_{*}^{\chi\neq 1}([X/G]).$$

Cette d\'ecomposition est fonctorielle pour les morphismes propres repr\'esentables.

De plus la projection sur le premier facteur correspond au morphisme naturel (\cite[Prop. $1.6$]{t2}) 
$\mathbf{G}_{*}([X/G])_{\mathbb{Q}} \longrightarrow \underline{\mathbf{G}}_{*}([X/G])_{\mathbb{Q}}$.

\end{cor}

\textit{Preuve:} On utilise \ref{t1}, et le fait que $C_{F}^{t} \longrightarrow F$ poss\`ede une section. 
Cette section donne un d\'ecomposition
$C_{F}^{t}\simeq F\coprod (C_{F}^{t})^{\neq 1}$, qui induit la d\'ecomposition voulue. Elle v\'erifie clairement les hypth\`eses
de fonctorialit\'e. 

La derni\`ere assertion provient du fait que $\psi_{F}$ commute avec la localisation \'etale (\ref{t2} $(4)$). $\Box$\\

On peut en r\'ealit\'e \^etre un peu plus explicite lorsque $S=Spec k$ est le spectre d'un corps, alg\'ebriquement clos 
pour simplifier, et $G$ est affine et lisse. \\

Notons $I_{F}^{t}$ le champ des ramifications mod\'er\'ees de $F$. C'est le champ classifiant les couples
$(s,h)$, o\`u $s$ est un objet de $F$ et $h$ est un automorphisme de $s$, d'odre premier \`a la caract\`eristique de $k$. 
On dispose d'un morphisme \'etale et fini $q : I_{F}^{t} \longrightarrow C_{F}^{t}$ qui envoie $(s,h)$ sur $(s,<h>)$, o\`u 
$<h>$ est le sous-groupe engendr\'e par $h$. On peut alors v\'erifier qu'il existe un isomorphisme de faisceau en $\overline{\mathbb{Q}}$-alg\`ebres sur $C_{F}^{t}$ ($\overline{\mathbb{Q}}$ est la clot\^ure alg\'ebrique de $\mathbb{Q}$)
$$q_{*}(\overline{\mathbb{Q}})\simeq \Lambda\otimes \overline{\mathbb{Q}}.$$
Cet isomorphisme n'est pas naturel car il demande de choisir un plongement $\mu_{\infty}(k) \hookrightarrow \mu_{\infty}(\overline{\mathbb{Q}})$,
que nous fixerons une bonne fois pour toute.

De l'isomorphisme ci-dessus, on d\'eduit un nouvel isomorphisme 
$$\underline{\mathbf{G}}_{*}(I_{F}^{t})_{\overline{\mathbb{Q}}}\simeq 
\underline{\mathbf{G}}_{*}^{\chi}(F)_{\overline{\mathbb{Q}}}.$$
En appliquant \ref{t1}, on trouve donc un isomorphisme 
$$\mathbf{G}_{*}(F)_{\overline{\mathbb{Q}}} \simeq \underline{\mathbf{G}}_{*}(I_{F}^{t})_{\overline{\mathbb{Q}}}.$$

Soit $T$ le sch\'ema des \'el\'ements de $G$ d'ordre fini et premier \`a la caract\'eristique de $k$. Par \cite[$XII$ $5.5$]{sga}, ce sch\'ema
s'\'ecrit comme une r\'eunion disjointe de ses orbites
$$T=\coprod_{g\in c(G)}G/Z_{g},$$
o\`u $c(G)$ est un ensemble de repr\'esentants de
l'ensemble des classes de conjugaisons d'\'el\'ements d'ordre fini et premier \`a $car k$, et 
$Z_{g}$ est le centralisateur de $g$ dans $G$. De plus, l'action de $G$ sur $T$ par conjugaison, correspond \`a l'action naturelle
de $G$ sur les espaces homog\`enes $G/Z_{g}$. On peut alors v\'erifier que $I_{F}^{t}$ est donn\'e par
$$I_{F}^{t}\simeq \coprod_{g\in c(G)}[X^{g}/Z_{g}]$$
o\`u $X^{g}$ est le sous-sch\'ema des points fixes de $g$ sur $X$.

On en conclut donc un isomorphisme 
$$\mathbf{G}_{*}(F)_{\overline{\mathbb{Q}}} \simeq \prod_{g \in c(G)}\underline{\mathbf{G}}_{*}([X^{g}/Z_{g}])_{\overline{\mathbb{Q}}}
\simeq \prod_{g \in c(G)}\mathbf{G}_{*}(X^{g}/Z_{g})_{\overline{\mathbb{Q}}}$$

\begin{cor}\label{c3}
Si $S=Spec k$ est le spectre d'un corps alg\'ebriquement clos, et $G$ est affine et lisse sur $k$, alors
il existe un isomorphisme
$$\mathbf{G}_{*}([X/G])_{\overline{\mathbb{Q}}} \simeq \prod_{g \in c(G)}\mathbf{G}_{*}(X^{g}/Z_{g})_{\overline{\mathbb{Q}}}.$$
\end{cor}

Rappelons aussi que l'on dispose du caract\`ere de Chern 
$$Ch : \mathbf{K}_{0}([X/G]) \longrightarrow  CH([X/G])_{\mathbf{Q}}:=H^{*}([X/G]_{et},\underline{K}_{*}\otimes \mathbb{Q}).$$
Ce caract\`ere de Chern se factorise par
$$\mathbf{K}_{0}([X/G]) \longrightarrow 
\underline{\mathbf{K}}_{0}([X/G])\longrightarrow  CH([X/G])_{\mathbf{Q}}.$$
Lorsque $X$ est de plus un sch\'ema r\'egulier, alors $Ch : \underline{\mathbf{K}}_{0}([X/H])_{\mathbf{Q}} \longrightarrow
CH([X/G])_{\mathbb{Q}}$ est un isomorphisme (\cite[Cor. $3.38$]{t2}). \\

Le corollaire suivant r\'epond positivement \`a la conjecture $(3.4)$ de \cite{v2}.

\begin{cor}\label{c4}
Si $X$ est r\'egulier, 
il existe un isomorphisme d'anneaux
$$\mathbf{K}_{0}([X/G])_{\mathbb{Q}}\simeq CH([X/G])_{\mathbb{Q}} \oplus 
\underline{\mathbf{K}}_{0}^{\chi\neq 1}([X/G]).$$
De plus la projection sur le premier facteur correspond au caract\`ere de Chern.
\end{cor}

\textit{Preuve:} Il suffit d'utiliser \ref{t1} et le fait que le carat\`ere de Chern induit un isomorphisme d'anneaux
$$\underline{\mathbf{K}}_{0}([X/H])_{\mathbf{Q}} \longrightarrow
CH([X/G])_{\mathbb{Q}}.$$
$\Box$ \\

Pour prendre en compte la $K$-th\'eorie sup\'erieure nous aurons besoin d'une petite digression sur les groupes de Chow sup\'erieurs
des champs alg\'ebriques. \\

Pour tout champ alg\'ebrique $F$ (toujours sur un corps), on dispose d'un complexe de pr\'efaisceau
$$\begin{array}{cccc}
\mathbf{Z}^{i} : & F_{et} & \longrightarrow & C(Ab) \\
& U & \mapsto & \mathbf{Z}^{i}(U)
\end{array}$$
o\`u $\mathbf{Z}^{i}(U)$ est le complexe de cycles de codimension $i$ d\'efini dans \cite{b}. On pose alors
$$CH^{i}(F,j):=H^{-j}(F_{et},\mathbf{Z}^{i}\otimes \mathbb{Q}).$$
On notera aussi $CH(F,\bullet):=\oplus_{i,j}CH^{i}(F,j)$. \\

Ces groupes de Chow sup\'erieurs poss\`edent de nombreuses bonnes propri\'et\'es. Nous les
\'enon\c{c}ons sans d\'emonstrations car le lecteur pourra trouver lui m\^eme les arguements sans aucunes difficult\'es. 
On doit aussi pouvoir trouver les arguments dans \cite{jo}.

\begin{enumerate}
\item Il existe des images directes par morphismes propres quelconques, ainsi que des images r\'eciproques pour
des morphismes plats quelconques. 

\item Les fonctorialit\'es ci-dessus satisfont aux formules usuelles de projection et de transfert. 

\item Il existe une suite exacte longue de localisation. 

\item Si $F$ est r\'egulier alors $CH(F,\bullet)$ devient un anneau bigradu\'e. De plus, il existe alors un caract\`ere de Chern
$$Ch : \underline{\mathbf{K}}_{m}(F) \longrightarrow CH(F,m)$$
qui rationnellement est un isomorphisme d'anneaux.

\item Pour tout champ alg\'ebrique, et $p : F \longrightarrow M$ sa projection sur son espace de modules, 
le morphisme $p_{*} : CH(F,\bullet) \longrightarrow CH(M,\bullet)$
est un isomorphisme. En particulier les groupes $CH(F,0)$ coincident modulo torsion avec les groupes 
de Chow usuel. 

\end{enumerate}

\textit{Remarque:} Nous avons d\'efini le caract\`ere de Chern, $Ch : \underline{\mathbf{K}}_{*}(F) \longrightarrow
CH(F,*)$. Cependant, nous pouvons composer avec le morphisme naturel $\mathbf{K}_{*}(F) \longrightarrow \underline{\mathbf{K}}_{*}(F)$, 
pour obtenir $Ch : \mathbf{K}_{*}(F) \longrightarrow CH(F,*)$. \\

Revenons au cas o\`u $F=[X/G]$ est le quotient d'un sch\'ema $X$ par un groupes affine et lisse. Alors \ref{t1} et 
les propri\'et\'es \'enonc\'ees ci-dessus impliquent le corollaire suivant.

\begin{cor}\label{c5}
Si $X$ est r\'egulier, 
il existe un isomorphisme d'anneaux gradu\'es
$$\mathbf{K}_{*}([X/G])_{\mathbb{Q}}\simeq CH([X/G],*) \oplus 
\underline{\mathbf{K}}_{*}^{\chi\neq 1}([X/G]).$$
De plus la projection sur le premier facteur correspond au caract\`ere de Chern.
\end{cor}

On peut m\^eme expliciter le terme $\underline{\mathbf{K}}_{*}^{\chi\neq 1}([X/G])$ en fonctions des groupes sup\'erieurs de Chow des points fixes. 
On obtient alors une g\'en\'eralisation 
de \cite[$3.4$]{v2} \`a la $K$-th\'eorie sup\'erieure.

\begin{cor}\label{c6}
Soit $F=[X/G]$ le champ quotient d'une sch\'ema en groupes affine et lisse sur un corps alg\'ebriquement clos, 
op\'erant sur un sch\'ema r\'egulier $X$. Notons $c(G)$ une ensemble de repr\'esentants de
l'ensemble des classes de conjugaisons d'\'el\'ements de $G$
d'ordre fini premier \`a la caract\'eristique de $k$. Pour $g \in c(G)$, notons aussi $X^{g}$ le sous-sch\'ema des points fixes, 
et $Z_{g}$ le centralisateur de $g$ dans $G$. Alors il existe un isomorphisme de $\overline{\mathbb{Q}}$-alg\`ebres
$$\mathbf{K}_{*}([X/G])_{\overline{\mathbb{Q}}}\simeq
\prod_{g \in c(G)}CH([X^{g}/Z_{g}],\bullet)_{\overline{\mathbb{Q}}}.$$
\end{cor} 

\textit{Preuve:} Rappelons qu'il existe un isomorphisme d'anneaux
$$\mathbf{K}_{*}(F)_{\overline{\mathbb{Q}}} \simeq \underline{\mathbf{K}}_{*}(I_{F}^{t})_{\overline{\mathbb{Q}}}.$$
En composant avec le caract\`ere de Chern $Ch : \underline{\mathbf{K}}_{*}(I_{F}^{t})_{\overline{\mathbb{Q}}} \simeq
CH(I_{F}^{t},*)_{\overline{\mathbb{Q}}}$, on en d\'eduit un isomorphisme
$$\mathbf{K}_{*}(F)_{\overline{\mathbb{Q}}}\simeq 
CH(I_{F}^{t},*)_{\overline{\mathbb{Q}}}.$$
On conclut alors en remarquant que $I_{F}^{t}\simeq \coprod_{g \in c(G)}[X^{g}/Z_{g}]$. $\Box$ \\

Enfin, lorsque $G$ est ab\'elien, le faisceau $\Lambda$ est alors un faisceau constant sur $C_{F}^{t}$. On peut alors en d\'eduire
facilement qu'il existe un isomorphisme de $\mathbb{Q}$-alg\`ebres (toujours pour $X$ r\'egulier)
$$\mathbf{K}_{*}(F)_{\mathbb{Q}}\simeq \prod_{c}CH([X^{c}/G],*)\otimes \mathbb{Q}(\zeta_{m(c)}),$$
o\`u le produit est pris sur l'ensemble des sous-groupes cycliques de $G$ et d'ordre premier \`a la caract\'eristique de $k$. 
Ceci r\'epond \`a la question $(3.6)$ de \cite{v2}. \\

Le corollaire \ref{c6} poss\`ede aussi une g\'en\'eralisation pour la $K$-th\'eorie rationnelle. La formule que l'on obtient est une
cons\'equence directe de \ref{t1}.

\begin{cor}\label{c6'}
Soit $F=[X/G]$ le champ quotient d'une sch\'ema en groupes affine et lisse sur un corps alg\'ebriquement clos, 
op\'erant sur un sch\'ema r\'egulier $X$. Notons $cycl(G)$ un ensemble de repr\'esentants de
l'ensemble des classes de conjugaisons des sous-groupes cycliques de $G$
d'ordre premier \`a la caract\'eristique de $k$. Pour $c \in cycl(G)$, notons aussi $X^{c}$ le sous-sch\'ema des points fixes, 
$N_{c}$ son normalisateur, $Z_{c}$ son centralisateur, et $W(c)=N_{c}/Z_{c}$ son groupe de Weyl.
Pour $c \in cycl(G)$, soit $\zeta_{m(c)}$ une racine de l'unit\'e d'ordre l'ordre de $c$.
Alors il existe un isomorphisme de $\mathbb{Q}$-alg\`ebres
$$\mathbf{K}_{*}([X/G])_{\mathbb{Q}}\simeq
\prod_{c \in cycl(G)}(CH([X^{c}/Z_{c}],\bullet)\otimes \mathbb{Q}(\zeta_{m(c)}))^{W(c)}.$$
\end{cor}

\textit{Preuve:} Il faut tout d'abord remarquer que 
$$C_{F}^{t}\simeq \coprod_{c \in cycl(G)}[X^{c}/N_{c}].$$
Ainsi, le th\'eor\`eme \ref{t1} donne un isomorphisme d'anneaux
$$\phi_{F} : \mathbf{K}_{*}([X/G])_{\mathbb{Q}}\simeq \prod_{c \in cycl(G)}\mathbb{H}^{-*}([X^{c}/N_{c}],\underline{\mathbf{K}}\otimes
\mathbb{Q}(\zeta_{m(c)})).$$
Comme $W(c)$ est un groupe fini, la formule de Leray appliqu\'ee au morphisme $[X^{c}/N_{c}] \longrightarrow [X^{c}/Z_{c}]$, donne
un isomorphisme naturel
$$\mathbb{H}^{-*}([X^{c}/N_{c}],\underline{\mathbf{K}}\otimes
\mathbb{Q}(\zeta_{m(c)})) \simeq (\mathbb{H}^{-*}([X^{c}/Z_{c}],\underline{\mathbf{K}}_{\mathbb{Q}})\otimes
\mathbb{Q}(\zeta_{m(c)})))^{W(c)}\simeq (\underline{\mathbf{K}}_{*}([X^{c}/Z_{c}])\otimes \mathbb{Q}(\zeta_{m(c)}))^{W(c)}.$$
On termine alors \`a l'aide de l'isomorphisme
$$Ch : \underline{\mathbf{K}}_{*}([X^{c}/Z_{c}]) \simeq CH([X^{c}/Z_{c}],*).$$
$\Box$

Supposons maintenant que $G$ soit r\'eductif, op\'erant sur un sch\'ema $X$ quelconque, toujours avec stabilisateurs r\'eduits et finis. On ne 
suppose plus n\'ecessairement que $k$ est alg\'ebriquement clos. 

La corollaire suivant r\'epond \`a la conjecture $(2.4)$ \cite{v2}. 
Notons qu'elle \'et\'e d\'emontr\'ee auparavant dans \cite[Cor. $5.2$]{eg}, et pour une situation un peu plus g\'en\'erale 
(sans hypoth\`eses sur $G$). 

\begin{cor}\label{c7}
Soit $f : [X/G] \longrightarrow BG$ le morphisme naturel.
Soit $\tau : \mathbf{G}_{0}([X/G]) \longrightarrow CH([X/G],0)$ la transformation de Riemann-Roch (\cite[$2.2$]{v2}). Alors 
$x \in Ker \tau$ si et seulement s'il existe $y \in \mathbf{K}_{0}(BG)$ de rang non-nul, tel que
$f^{*}(y).x=0$.
\end{cor}

\textit{Preuve:} Par la propri\'et\'e du module (\ref{t2} $(5)$), la condition est \'evidemment suffisante. \\

Par une descente galoisienne, on peut clairement supposer que $k$ est alg\'ebriquement clos. 

On commence par remarquer que le noyau de $\tau$ est aussi le noyau du morphisme canonique
$\mathbf{G}_{0}([X/G]) \longrightarrow \underline{\mathbf{G}}_{0}([X/G])$. En effet, on peut v\'erifier que le morphisme $\tau$
se factorise par (car le caract\`ere de Chern se factorise de cette fa\c{c}on)
$$\xymatrix{
\mathbf{G}_{0}([X/G]) \ar[r]^-{\tau} \ar[d]_-{can} & CH([X/G])  \\
\underline{\mathbf{G}}_{0}([X/G]) \ar[ru]_-{\underline{\tau}} & }$$
On conclut alors en remarquant que $\underline{\tau}$ est un isomorphisme (\cite[$3.38$]{t2}). 

Par \ref{c2}, on voit donc que 
$x \in \underline{\mathbf{G}}_{0}^{\chi\neq 1}([X/G])$. En utilisant la propi\'et\'e du module (\ref{t2} $(5)$), on
voit qu'il suffit  de trouver 
un \'el\'ement $y \in \mathbf{K}_{0}(BG)$, tel que $\phi_{F}(f^{*}y) \in \underline{\mathbf{K}}_{0}([X/G])
\hookrightarrow \underline{\mathbf{K}}_{0}^{\chi}([X/G])$. 
En analysant la construction de $\phi_{[X/G]}$, 
il est facile de v\'erifier qu'il est suffisant de trouver $y \in \mathbf{K}_{0}(BG)$ de rang non-nul, tel que son caract\`ere
$$\chi(y) : G \longrightarrow k$$
soit nul sur tout \'el\'ement non trivial de $G$ qui fixe au moins un point de $X$, et dont l'ordre est premier \`a la caract\'eristique.

Soit $S$ l'ensemble des tels \'el\'ement. Alors $G$ op\`ere sur $S$ par conjugaison, et l'ensemble quotient est un ensemble fini. 
Choisissons un syst\`eme de rep\'esentant $S'$ dans $G$. Si l'on construit pour tout $h \in S'$, un \'el\'ement $y_{h} \in \mathbf{K}_{0}(BG)$
de rang non-nul, et dont le caract\`ere s'annule en $h$, l'\'el\'ement $y=\prod_{h \in S'}h_{y}$ r\'epondra \`a la question. 

Soit $h \in S'$. Comme $h$ est d'ordre fini premier \`a $car k$, il est semi-simple, et est donc contenu dans un tore
maximal $h \in \mathbb{T}$. Le morphisme induit sur les groupes des caract\`ers $\mathcal{X}(\mathbb{T}) \longrightarrow
\mathcal{X}(<h>)$ est surjectif. Il existe donc $z \in \mathbf{K}_{0}(B\mathbb{T})$, de rang $1$, et 
dont la restriction \`a $<h>$ est une repr\'esentation fid\`ele. On consid\`ere $1+z+\dots+z^{m-1}$, o\`u $m$ est l'ordre
de $h$. C'est un \'el\'ement de rang non-nul de $\mathbf{K}_{0}(B\mathbb{T})$, dont le caract\`ere s'annule en $h$. Notons le encore $z$.
Notons $W$ le groupe de Weyl de $\mathbb{T}$. Alors, en remplac\c{c}ant
$z$ par $\sum_{\sigma \in W}\sigma(z)$, on peut supposer que $z$ est invariant par le groupe de Weyl. Il provient donc
d'un $y_{h} \in \mathbf{K}_{0}(BG)_{\mathbb{Q}}$, dont un multiple est l'\'el\'ement cherch\'e. $\Box$\\

Lorsque $F$ n'est plus un quotient, on peut toujours se poser la question de la validit\'e des corollaires
\ref{c2}, \ref{c5} et \ref{c7}. Les deux premier se g\'en\'eralisent sans aucuns probl\`emes. Pour ce qui est de
\ref{c7}, on peut r\'epondre tout au moins dans le cas r\'egulier. 

\begin{cor}\label{c8}
Soit $F$ un champ alg\'ebrique r\'egulier (plus n\'ecessairement sur un corps). Notons $can : \mathbf{G}_{*}(F)_{\mathbb{Q}} \longrightarrow
\underline{\mathbf{G}}_{*}(F)$ le morphisme naturel. Alors, $x \in Ker(can)$ si et seulement s'il existe
$y \in \mathbf{K}_{0}(F)$ de rang non-nul, tel que $y.x=0$. 
\end{cor} 

\textit{Preuve:} De nouveau la propri\'et\'e du module implique que la condition est suffisante. \\

Par \ref{t1} le noyau de $can$ est $\underline{\mathbf{G}}_{*}^{\chi\neq 1}(F)$, et par \ref{t2} il nous suffit donc de prendre 
l'image inverse par $\phi_{F}$ de n'importe quel \'el\'ement de rang non nul dans $\underline{\mathbf{K}}_{0}(F) \hookrightarrow 
\underline{\mathbf{G}}_{0}^{\chi}(F)$. $\Box$\\

\textit{Remarque:} On peut aussi d\'emontrer le corollaire \ref{c8} de fa\c{c}on tout \`a fait \'el\'ementaire
en utilisant le m\^eme arguement que celui utilis\'e pour la preuve de \cite[Thm. $3.1$]{v2}. Cependant, le d\'emonstration
que l'on en donne permet de voir qu'il suffit de savoir r\'epondre \`a la questions suivante. 

Soit $F$ un champ alg\'ebrique. Existe-t-il un \'el\'ement $x \in \mathbf{K}_{0}(F)$ tel que 
$\phi_{F}(x) \in \underline{\mathbf{K}}_{0}(F) \hookrightarrow \underline{\mathbf{K}}^{\chi}_{0}(F)$, avec
$\phi_{F}(x) \in \underline{\mathbf{K}}_{0}(F)$ de rang non-nul (par exemple $\phi_{F}(x)=1$) ? \\

Notons que pour y r\'epondre il faudra certainement savoir
r\'epondre \`a la question suivante. C'est une g\'en\'eralisation de la propri\'et\'e de descente de la $K$-th\'eorie
rationnelle d\'emontr\'ee par Thomason dans \cite[Thm. $11.10$]{th}. \\

\textit{Question:} Soit $p : F \longrightarrow M$ la projection d'un champ alg\'ebrique sur son espace de modules.
Le morphisme naturel
$$\mathbf{K}_{*}(F)_{\mathbb{Q}} \longrightarrow \mathbb{H}^{-*}(M_{et},p_{*}\mathbf{K}_{\mathbb{Q}})$$
est-il un isomorphisme ? \\

Il me semble que d\'ej\`a le cas o\`u $F$ est un espace alg\'ebrique n'est pas clair. \\

\textit{Remarque:} Lorsque $S=Spec k$, avec $k$ un corps de caract\'eristique nulle, il existe une question analogue pour
la cohomologie p\'eriodique des champs alg\'ebrique (voir la remarque qui suit \cite[Prop. $3.54$]{t2}). \\

Il me semble que la description des anneaux $\mathbf{K}_{*}([X/H])_{\mathbf{Q}}$, o\`u $H$ est un groupe fini op\'erant sur un sch\'ema $X$, est 
aussi li\'ee \`a ces questions, et c'est sans doute le premier probl\`eme auquel il est important d'apporter une r\'eponse. 

\end{subsection}

\begin{subsection}{Nouvelles $\gamma$-filtrations et formalisme de Riemann-Roch}

Soit $F$ un champ alg\'ebrique, et $\mu$ le faisceau sur $C_{F}^{t}$ des racines de l'unit\'e d'ordre correspondant au sous-groupe
cyclique. En clair, si $U \longrightarrow C_{F}^{t}$ est un morphisme \'etale, correspondant \`a un objet $s \in F(U)$ et 
$c$ un sous-groupe cyclique de $\underline{Aut}_{U}(s)$, $\mu(U):=\mu_{m(c)}(U)$, o\`u $m(c)$ est l'ordre de $c$. Comme 
$m(c)$ est inversible sur $U$, $\mu$ est un sch\'ema en groupes \'etale et fini sur $C_{F}^{t}$. 
C'est m\^eme un sch\'ema en groupes cyclique de type multiplicatifs. Notons $\mathbb{Q}[\mu]$ le faisceau en $\mathbb{Q}$-alg\`ebres 
de groupes, et $\mathbb{Q}[\mu] \longrightarrow \Gamma$ le quotient correspondant localement au morphisme naturel
$\mathbb{Q}[\mathbb{Z}/m] \longrightarrow \mathbb{Q}(\zeta_{m})$. \\

Notons $q : I_{F}^{t} \longrightarrow C_{F}^{t}$ le morphisme naturel. Alors on peut v\'erifier qu'il existe un isomorphisme
de faisceaux en $\mathbb{Q}$-alg\`ebres
$$\Lambda\otimes \Gamma \simeq q_{*}q^{*}(\Gamma).$$
Ceci implique qu'il existe un isomorphisme
$$\mathbb{H}((C_{F}^{t})_{et},\underline{\mathbf{K}}_{\Lambda\otimes \Gamma}) \simeq
\mathbb{H}((I_{F}^{t})_{et},\underline{\mathbf{K}}_{q^{*}\Gamma}).$$
Soit $S' \longrightarrow S$ un rev\^etement galoisien ramifi\'e de groupe $H$, tel que le faisceau $\mu$ soit constant 
sur $C_{F}^{t}\times_{S}S'$. On peut prendre par exemple $S'=\mu_{m}$, o\`u $m$ est un multiple de tous les ordres des groupes d'isotropie
de $F$. Sur une composante connexe de $C_{F}^{t}$ sur laquelle le groupe cyclique universel est d'ordre $m$, le faisceau
$\Gamma$ devient alors isomorphe au faisceau constant $\mathbb{Q}(\zeta_{m})$. 

Ceci implique qu'il existe un isomorphisme de $\mathbb{Q}$-alg\`ebres
$$\mathbb{H}^{-*}(C_{F}^{t}\times_{S}S',\underline{\mathbf{K}}_{\Lambda\otimes \Gamma}) \simeq
\underline{\mathbf{K}}_{*}(I_{F}^{t}\times_{S}S')\otimes \Gamma.$$
De plus cet isomorphisme est \'equivariant pour l'action de  $H$ sur $S'$ et $\Gamma$. En fin de compte on trouve un morphisme naturel
de $\mathbb{Q}$-alg\`ebres
$$\hspace*{-15mm}
\xymatrix{ \mathbf{K}_{*}(F) \ar[r]^-{\phi_{F}} & \underline{\mathbf{K}}_{*}^{\chi}(F) \ar[r] &
\underline{\mathbf{K}}_{*}^{\chi}(F\times_{S}S')^{H} \ar[r] & \mathbb{H}^{-*}(C_{F}^{t}\times_{S}S',\underline{\mathbf{K}}_{\Lambda\otimes \Gamma})^{H} \ar[r] & (\underline{\mathbf{K}}_{*}(I_{F}^{t}\times_{S}S')\otimes \Gamma)^{H}}$$
En utilisant \cite[Thm. $3$]{gs}, $\underline{\mathbf{K}}_{*}(I_{F}^{t}\times_{S}S')$ est un $\lambda$-anneau, 
et on peut trouver un isomorphisme d'anneaux (car $*$ est $\underline{\mathbf{K}}_{\mathbb{Q}}$-coh\'erent sur $F_{et}$, \cite[Prop. $8$]{gs})
$$\underline{\mathbf{K}}_{*}(I_{F}^{t}\times_{S}S') \simeq Gr_{\gamma}^{.}\underline{\mathbf{K}}_{*}(I_{F}^{t}\times_{S}S').$$
Cet isomorphisme \'etant fonctoriel, il est $H$-\'equivariant, et permet donc de trouver un isomorphisme naturel de $\mathbb{Q}$-alg\`ebres
$$(\underline{\mathbf{K}}_{*}(I_{F}^{t}\times_{S}S')\otimes \Gamma)^{H} \simeq
(Gr_{\gamma}^{.}\underline{\mathbf{K}}_{*}(I_{F}^{t}\times_{S}S')\otimes \Gamma)^{H}.$$
Remarquons que l'action de $H$ sur le second membre est compatible avec la graduation.

\begin{df}\label{d6}
Pour tout champ alg\'ebrique $F$, nous noterons
$$H^{p}_{\chi}(F,j):=(Gr_{\gamma}^{j}\underline{\mathbf{K}}_{2j-p}(I_{F}^{t}\times_{S}S')\otimes \Gamma)^{H}.$$
La $\mathbb{Q}$-alg\`ebre $\bigoplus_{p,j}H^{p}_{\chi}(F,j)$ sera appel\'ee l'alg\`ebre de cohomologie motivique de $F$ \`a coefficients
dans les caract\`eres.
\end{df}

Notons que par construction, il existe un morphisme d'anneaux 
$$Ch^{\chi} : \mathbf{K}_{*}(F) \longrightarrow H^{*}_{\chi}(F,\bullet).$$

Supposons maintenant que $F$ soit quasi-projectif. Notons $MI_{F}^{t}$ l'espace de modules de $I_{F}^{t}$. Alors, par les m\^emes arguments
que ci-dessus, ainsi que par \ref{t1}, on peut trouver un isomorphisme fonctoriel pour les images directes
$$\mathbf{G}_{*}(F) \longrightarrow (Gr^{\gamma}_{.}\mathbf{G}_{*}(MI_{F}^{t}\times_{S}S')\otimes \Gamma)^{H},$$
o\`u $Gr^{\gamma}_{.}$ est la $\gamma$-filtration d\'efinie dans \cite{so} (qui existe car $MI_{F}^{t}$ est quasi-projectif par
hypoth\`ese).

\begin{df}\label{d7}
Pour tout champ alg\'ebrique quasi-projectif $F$, nous noterons
$$H_{p}^{\chi}(F,j):=(Gr_{j}^{\gamma}\mathbf{G}_{p-2j}(MI_{F}^{t})\otimes \Gamma)^{H}.$$
Le groupe $\bigoplus_{p,j}H^{\chi}_{p}(F,j)$ sera appel\'e le groupe d'homologie motivique de $F$ \`a coefficients dans les
repr\'esentations.
\end{df}

Aussi par construction, il existe un morphisme covariant pour les images directes de morphismes propres et repr\'esentables (ou 
de morphismes propres de dimension cohomologique fini quelconque si $F$ est lisse)
$$\tau_{F}^{\chi} : \mathbf{G}_{*}(F) \longrightarrow H^{\chi}_{*}(F,\bullet).$$

Ces deux morphismes
$$Ch^{\chi} : \mathbf{K}_{*}(F) \longrightarrow H^{*}_{\chi}(F,\bullet)$$
$$\tau_{F}^{\chi} : \mathbf{G}_{*}(F) \longrightarrow H_{*}^{\chi}(F,\bullet)$$
d\'efinissent par images r\'eciproques des filtrations sur les groupes $\mathbf{K}_{m}(F)$ et $\mathbf{G}_{m}(F)$. Le lecteur
se convaincra (en consid\'erant l'exemple de $F=[X/H]$ par exemple) que ces filtrations ne sont pas \'equivalentes (m\^eme rationnellement)
aux $\gamma$-filtrations 'usuelles' (i.e. provenant des puissances ext\'erieures de fibr\'es sur $F$). De m\^eme, la filtration
induite sur $\mathbf{G}_{0}(F)$ n'est pas \'equivalente \`a la filtration par la dimension du support. La raison en est que
$C_{F}^{t}$ n'est pas \'equidimensionel en g\'en\'eral. 

Bien entendu la construction que nous donnons de ces filtrations n'est pas satisfaisante, dans le sens o\`u il devrait exister un moyen 
de le d\'efinir plus directement sur les groupes de $K$-th\'eorie, sans passer par les morphismes $\phi$ et $\psi$. \\

\textit{Question:} Existe-t-il des descriptions de ces filtrations intrens\`eques \`a $F$ ? \\

Supposons maintenant que $F$ soit lisse sur $S$. On a d\'ej\`a vu qu'il existe alors un \'el\'ement inversible
$\alpha_{F} \in \underline{\mathbf{K}}_{0}^{\chi}(F)$.
L'image de $\alpha_{F}$ par le morphisme naturel
$$\underline{\mathbf{K}}_{0}(F) \longrightarrow H^{*}_{\chi}(F,2*)$$
sera not\'ee $Ch(\alpha_{F})$. Comme $\alpha_{F}$ est  un \'el\'ement
inversible, $Ch(\alpha_{F})$ aussi.

Le champ $C_{F}^{t}$ \'etant encore lisse sur $S$, il poss\`ede un fibr\'e tangent $T$, qui poss\`ede un classe dans
$\underline{\mathbf{K}}_{0}(C_{F}^{t})$. Comme $\underline{\mathbf{K}}_{0}(C_{F}^{t})$ est un $\lambda$-anneau, on dispose
de l'application 'classe de Todd'
$$Td : \underline{\mathbf{K}}_{0}(C_{F}^{t}) \longrightarrow 
Gr_{\gamma}^{*}\underline{\mathbf{K}}_{0}(C_{F}^{t}).$$
L'image de $Td(T)$ par le morphisme naturel $Gr_{\gamma}^{*}\underline{\mathbf{K}}_{0}(C_{F}^{t}) \longrightarrow
H^{*}_{\gamma}(F,2*)$ sera not\'ee $Td(C_{F}^{t})$.

\begin{df}\label{d8}
Si $F$ est un champ alg\'ebrique lisse, alors sa classe de Todd est d\'efinie par
$$Td^{\chi}(F):=Ch(\alpha_{F})^{-1}.Td(C_{F}^{t}) \in H^{*}_{\chi}(F,2*).$$
\end{df}

En utilisant alors les techniques de d\'emonstration du th\'eor\`eme \ref{t2}, 
on peut esp\'erer d\'emontrer un th\'eor\`eme de Riemann-Roch \`a valeurs dans cette cohomologie
motivique. Comme nous n'avons pas v\'erifi\'e les d\'etails, nous laisserons cette question ouverte. \\

\textit{Probl\`emes:} \begin{enumerate}
\item Montrer que $F \mapsto (H^{*}_{\chi}(F,\bullet),H^{\chi}_{*}(F,\bullet))$ forme une 'bonne' (cf ci-dessous) th\'eorie
cohomologique. 
\item Montrer alors que les th\'eor\`emes \cite[$4.11$]{t1} et \cite[$3.33$]{t2} 
se g\'en\'eralisent au cas des champs alg\'ebriques quasi-projectifs sur $S$. 
\end{enumerate}

\textit{Remarque:} Dans le premier probl\`eme, il me semble peu raisonable d'attendre que la th\'eorie satisfasse \`a tous
les axiomes de \cite{so}. Ceci est du essentiellement au fait que la th\'eorie n'existe pas \`a coefficients dans $\mathbb{Q}$,
mais uniquement apr\`es l'extension des scalaires \`a $\Gamma$. Ainsi, il me semble difficile de poss\`eder d'une bonne th\'eorie des classes
de Chern. Ceci est aussi \`a raprocher du fait que lorsque $F$ est connexe, $C_{F}^{t}$ ne l'est pas, et donc on perd la notion
m\^eme de rang. 

\end{subsection}

\end{section}

\begin{section}{Cas des champs d'Artin}

La plupart des constructions que nous avons donn\'ees dans les paragraphes pr\'ec\'edents gardent en r\'ealit\'e un sens lorsque l'on remplace
les champs de Deligne-Mumford par des champs d'Artin. Ainsi, pour tout champ d'Artin (de type fini sur $S$), on dispose du champ $C_{F}^{t}$
classifiant les paires $(s,c)$, o\`u $s$ est un objet de $F$ au-dessus d'un sch\'ema $X$, et $c$ est un sous-groupe cyclique
de type multiplicatif de l'espace alg\'ebrique en groupes $\underline{Aut}_{X}(s)$. En mettant de cot\'e les cas pathologiques, $C_{F}^{t}$ est
cetainement encore un champ d'Artin, localement de type fini sur $S$. Sur $C_{F}^{t}$ on dispose toujours du champ en groupes universel, et 
de son faisceau des caract\`eres $\mathcal{X}$. On peut donc donner un sens \`a $\Lambda$. Le lecteur v\'erifiera sans probl\`emes que
la construction donn\'ee dans le pr\'esent travail se g\'en\'eralise, et donne un morphisme
$$\phi_{F}  : \mathbf{K}(F) \longrightarrow \mathbb{H}(C_{F}^{t},\underline{\mathbf{K}}_{\Lambda}).$$
Le simple exemple du champ classifiant d'un groupe alg\'ebrique montre que \ref{t1} $(1)$ n'a aucune chance de rester vrai sans hypoth\`eses
suppl\'ementaires. Il semble par exemple raisonnable de commencer par se restreindre au cas des champs dont la diagonale
est un morphisme fini. Le r\'esultat d\'emontr\'e dans \cite{vv} montre qu'il est alors tr\`es probable 
que \ref{t1} $(1)$ reste vrai, mais je ne sais pas le d\'emontrer. 

Je ne vois pas comment g\'en\'eraliser la construction de $\psi_{F}$. \\

En ce qui concerne le th\'eor\`eme \ref{t2}, ou encore une \'eventuelle formule de Riemann-Roch, 
on peut raisonnablement attendre \`a ce qu'il soit vrai pour un morphisme
propre de dimension cohomologique fini entre deux champs dont les diagonales sont affines. En effet, dans ce cas on se trouve avec une
famille de champs dont la diagonale est finie, parametr\'ee par un champ d'Artin. Pour les morphismes fortement projectif, les
techniques standards donnent une formule de Riemann-Roch (\cite[$3.23$]{t2}).

\end{section}

\end{document}